\documentclass[a4paper,12pt]{article}
\usepackage{amsmath}
\usepackage{amssymb, amsbsy, amstext}
\usepackage{srcltx}
\usepackage{amscd,amsfonts,color}
\usepackage[figuresright]{rotating}
\usepackage{url}

\input amssym.def
\input amssym.tex
\newtheorem{lem}{Lemma}[section]%
\newtheorem{theorem}[lem]{Theorem}%
\newtheorem{defi}[lem]{Definition}%
\newtheorem{cor}[lem]{Corollary}%
\newtheorem{prop}[lem]{Proposition}%
\newtheorem{rem}[lem]{Remark}%
\newtheorem{conj}[lem]{Conjecture}%

\def\a{\alpha} \def\b{\beta} \def\g{\gamma} \def\d{\delta} \def\e{\varepsilon}

\def\si{\Sigma} \def\O{\Omega} \def\G{\Gamma}
\def\D{{\rm D}}

\def\di{\bigm|} \def\lg{\langle} \def\rg{\rangle}

\def\PSL{\hbox{\rm PSL}}  
\def\Aut{\hbox{\rm Aut\,}} \def\Inn{\hbox{\rm Inn}} \def\Syl{\hbox{\rm Syl}}
 \def\soc{\hbox{\rm soc}} \def\Fix{\hbox{\rm Fix }}

 \def\PG{\hbox{\rm PG}} \def\PGL{\hbox{\rm PGL}}\def\PSU{\hbox{\rm PSU}}
  \def\GL{\hbox{\rm GL}}  \def\P\GL{\hbox{\rm P\GL}}
 \def\SL{\hbox{\rm SL}} \def\FF{{\hbox{\sf F\kern-.43emF}}}
  
\def\PGammaL{{\rm P\Gamma L}} 
 \def\Inndiag {{\rm Inndiag}}
\def\Sym{\hbox{\rm Sym}}

\def\N{\hbox{\rm N}}
\def\C{\hbox{\rm C}}
\def\Z{\hbox{\rm Z}}
\def\A{\hbox{\rm A}}
\def\soc{\hbox{\rm soc}}

\def\Sz{\hbox{\rm Sz}}
\def\Ree{\hbox{\rm Ree}}

\def\Out{\hbox{\rm Out}}

\def\ZZ{\mathbb{Z}} 

\def\nd{\mathrel{\bigm|\kern-.7em/}} 
 \def\f{\noindent}
\def\qed{\hfill $\Box$} \def\demo{\f {\bf Proof}\hskip10pt}

\begin{document}
 \begin{center} {\bf\large The Burness-Giudici Conjecture on Primitive Groups\\ with Socle $\Ree(q)$ and $\Sz(q)$}
\vskip 3mm
{\sc Huye Chen}\\
{\footnotesize
School of Mathematics \& Center for Applied Mathematics of Guangxi,\\ Guangxi University,   Nanning 530004, P. R. China}\\

{\sc Shaofei Du}\\
{\footnotesize
School of Mathematical Sciences, Capital Normal University, Beijing 100048, P. R. China}\\

\end{center}

\renewcommand{\thefootnote}{\empty}
 \footnotetext{{\bf Keywords} Saxl graph, base of permutation group,  simple group}
  \footnotetext{E-mail addresses:  chenhy280@gxu.edu.cn (H. Chen), dushf@mail.cnu.edu.cn (S. Du).}

\begin{abstract}
Let $G$ be a transitive permutation group on $\O$ containing two points $\a, \b$ such that
$G_\a\cap G_\b=1$. The  Saxl graph $\Sigma(G)$ of $(G, \O)$ is defined as the graph with vertex set $\O$, where two vertices $\a', \b'$ are adjacent if and only if $G_{\a'}\cap G_{\b'}=1$.
Burness and Giudici conjectured that for any primitive permutation group $G$, its Saxl graph $\Sigma(G)$ satisfies the property that any two vertices share a common neighbor.
We study this common-neighbor property for primitive groups whose socle is a simple group of Lie type of rank one, namely $\PSL(2,q)$, $\PSU(3,q)$, $\Ree(q)$ or $\Sz(q)$. In this paper, we establish the property for groups with socle $\Ree(q)$ or $\Sz(q)$, except possibly when $\soc(G)=\Ree(q)\lneq G $ and a point stabilizer $M$ satisfies $M\cap \soc(G)\cong\ZZ_2\times\PSL(2,q)$.
       \end{abstract}

\section{Introduction}
Let $G$ be a permutation group on a finite set $\Omega$. A \emph{base} for $G$ is a subset of $\Omega$ such that its pointwise stabilizer is trivial.  The \emph{base size} $b(G)$ of $G$  is the minimal
size of a base for $G$, while a \emph{base size set} $\Delta$ is just a base with  size $b(G)$.
In \cite{BG} Burness and Giudici  introduced  the  \emph{Saxl graph} $\Sigma(G)$ of  a base-two permutation group $G$ on $\O$,  with the vertex set
$\Omega$, while  two vertices are adjacent if and only if they form a base for $G$.
 From now on, we assume that $G\leq \Sym(\O)$ acts transitively on $\O$. Fixing a point $\alpha\in \Omega$, we call the orbits of $G_\a$ the {\em suborbits} of $G$ relative to $\a$; the trivial suborbit is $\{\a\}$ itself.
   Now $\{\alpha,\beta\}$ is a base for $G$ with $b(G)=2$ if and only if $G_{\alpha}$ acts regularly on the suborbit containing $\b$.
   It follows that the neighborhood $\si_1(\a)$ of $\alpha$ in
  $\Sigma(G)$ is the union $\G$ of all regular suborbits of $G$ relative to $\a$. Burness and Giudici \cite{BG} conjectured that
\begin{conj}$($BG-Conjecture$)$ \label{mainC'}
If $G$ is a finite primitive permutation group with $b(G)=2$, then every pair of vertices in its Saxl graph $\Sigma(G)$ has a common neighbor.
\end{conj}
They used probabilistic methods to verify it for a broad class of primitive groups. Specifically, their work establishes the conjecture for symmetric and alternating groups $G = S_n$ or $A_n$ (with $n > 12$) whose point stabilizer is primitive, as well as for certain primitive groups of diagonal type and sufficiently large twisted wreath products. These latter cases were later extended by Huang (see Sections 5.8 and 5.9 of \cite{HPHD}).
Furthermore, using the computational algebra system Magma, they verified the conjecture for all primitive groups of degree at most $4095$. The same paper \cite{BG} also confirms the conjecture for many sporadic simple groups.
Subsequently, Burness and Huang \cite{BH} proved the conjecture for almost simple primitive groups with soluble point stabilizers. In a separate paper \cite{BH2}, the same authors further investigated the conjecture for primitive groups of product type.
More recently, Lee and Popiel \cite{LP} provided a proof for most affine-type groups with sporadic point stabilizers. The framework of the conjecture has since been expanded by the introduction of generalized Saxl graphs and Saxl hypergraphs in \cite{FH, LP1}.

However, Rizzoli and Thomas \cite{RT} have recently constructed counterexamples to the BG-Conjecture of affine, product and twisted wreath type. Thus, the conjecture is false in general, although its restriction to almost simple primitive groups remains open.

To investigate the BG-Conjecture for almost simple primitive groups whose socle is a finite simple group of Lie type, it is natural to begin with the rank-one families.
 All such groups include $\PSL(2, q)$, $\Ree(q)$, $\Sz(q)$ and
$\PSU(3,q)$. The case in which the socle is  $\PSL(2, q)$ was treated by the present authors  ~\cite{Chen-Du} when $\{\a, \b\}$ is not a base, and by Burness and Huang~\cite{BH} when  $\{\a, \b\}$ is a base.
In this paper, we consider the case in which $\soc(G)\in \{\Ree(q),\Sz(q)\}$.
Our main result is the following.
\begin{theorem}\label{main}
Let $G$ be a finite primitive permutation group with base size $2$, and suppose that its socle is either $\Sz(q)$ or $\Ree(q)$. Then the BG-Conjecture holds for the Saxl graph $\Sigma(G)$, except possibly when  $\soc(G)=\Ree(q)\lneq G$ and a point stabilizer $M$ satisfies $M\cap \soc(G)\cong\ZZ_2\times \PSL(2,q).$
\end{theorem}
We also establish the following base-size result for primitive groups with insoluble point stabilizers.
\begin{theorem}
Let $G$ be a finite almost simple primitive permutation group with socle $\Ree(q)$ or $\Sz(q)$, and let $M$ be a point stabilizer. If $M$ is insoluble, then $b(G)=2.$
\end{theorem}
This follows from the maximal subgroup classifications, Theorem 3.8, Lemma 3.17 and Theorem 4.9, together with Remark 2.4.
\section{Preliminary Results}
For a subgroup $H$ of a group $G$, we write $\N_G(H)$ and $\C_G(H)$ for the normalizer and centralizer of $H$ in $G$, respectively. As usual, $\Z(G):=\C_G(G)$.
Let $\Omega:=[G:H]$
be the set of right
cosets of $H$ in $G$, on which $G$ acts by right multiplication. For $K\leq G$, define  $\Fix_{\Omega}(K):=|\{\omega\in\Omega\mid \omega^k=\omega~~\text{for every}~~k\in K\}|.$ When the underlying set $\Omega$ is clear, we simply write $\Fix(K)$. For $\alpha\in\Omega$, we write $G_{\alpha}:=\{g\in G\mid \alpha^g=\alpha\}$ for the stabilizer of $\alpha$ in $G$.
Throughout this paper, a diagonal matrix with diagonal entries $a_1,a_2,\cdots,a_n$ is denoted by $[a_1,a_2,\cdots,a_n];$ and an anti-diagonal matrix whose anti-diagonal entries, read from the lower-left corner to the upper-right corner, are $b_1,b_2,\cdots,b_n$ is denoted by $]b_1,b_2,\cdots,b_n[$. For a finite group $G$ and a prime $r$, we denote by $\Syl_r(G)$ the set of Sylow $r$-subgroups of $G$.

\begin{prop}\label{man}
{\rm \cite{Man}}. \
Let $G$ be a finite group acting transitively on $\O$,  and let $H:=G_\a$
for some $\a\in \O$.
Suppose that $K\le G$ has a $G$-conjugate contained in $H$. Assume that the $G$-conjugates of $K$ contained in $H$ form $t$ conjugacy classes under $H$,  with representatives $K_1$, $K_2$, $\cdots,$ $K_t$.
Then $K$ fixes $\sum_{i=1}^{t}|\N_G(K_i):\N_H(K_i)|$ points of $\O$.
\end{prop}


For a simple group $T$ of Lie type, we define $\Inndiag(T)$ to be the group of inner-diagonal automorphisms of $T$, that is, the subgroup of $\Aut(T)$ generated by the inner and diagonal automorphisms.
The following proposition is derived from \cite[p.79]{BGbook}.
\begin{prop}\label{fieldconjugate}
Let $T$ be a simple group of Lie type and let $x\in \Aut(T)$ be a field or graph-field automorphism of prime order $r$. Then $y\in \Inndiag(T)x$ has order $r$ if and only if $x$ and $y$ are
$\Inndiag(T)$-conjugate.
\end{prop}

The following proposition is a modification of
results from \cite[Lemmas 3.6 and 3.8]{BG}.

\begin{prop}\label{prob}
Let $G$ be a finite transitive nonregular permutation group acting by right multiplication on the coset space $\O:=[G:M]$. Let $\check{Q}(G):=\frac{|M|}{|\O|}(\sum_{i=1}^k(|\lg x_i\rg|-1)\frac{\Fix(\lg x_i\rg)}{|\N_M(\lg x_i\rg )|})$, where
 $\mathcal{P}^*(M):=\{\lg x_1\rg,\lg x_2\rg,\cdots,\lg x_k\rg\}$ is a set of representatives for the $M$-conjugacy classes of subgroups of prime order in $M$ and $\Fix(\lg x_i\rg)$ is the number of
 points of $\O$ fixed by $\lg x_i\rg$.  Then $b(G)=2$ and any two vertices in $\Sigma(G)$ have a common neighbor, provided  $\check{Q}(G)< \frac 12$.
\end{prop}

\begin{rem}\label{big}
 Let $G_0\leq G\leq \Sym(\Omega)$ be transitive nonregular permutation groups. If $b(G)=2$, then $b(G_0)=2$ and $\Sigma(G)$ is a spanning subgraph of $\Sigma(G_0)$. Consequently, if every two vertices of $\Sigma(G)$ have a common neighbor, then the same conclusion holds for $\Sigma(G_0)$.
\end{rem}

\section{Groups with $\soc(G)=\Ree(q)$}
To prove Theorem~\ref{main} in the case where $T:=\soc(G)=\Ree(q)$, we need to consider all primitive permutation representations of $G$ on $[G:M]$, where $M$ ranges over the core-free maximal subgroups of $G$.
  In \cite{BH}, Burness and Huang proved the conjecture in the case  when $M$ is  soluble.   Therefore, we assume that $M$ is insoluble.
    Checking the structures of $\Ree(q)$ given in \cite[Theorem $4.2$]{Wilson},  the remaining possibilities are $M\cap \soc(G)=\ZZ_2\times\PSL(2,q)$ or $\Ree(q'),$  where $q=3^m$ and $q'=3^{m'}$, $m=m'e$ and $e$ is an odd prime.
Recall that    $\Out(\soc(G))=\lg f\rg $, where $f$ is a field automorphism.
In Subsection 3.1, some preliminaries  for $\Ree(q)$ will be given and in Subsections 3.2 and 3.3, two possibilities for $M_0:=M\cap T$ will be dealt with separately.
Unless otherwise stated, all notation and standing assumptions introduced above or in Subsection 3.1 will remain in force throughout this section.

\subsection{Preliminaries for $\Ree(q)$}\label{Ree}
By \cite[Lemma 3]{LN},  $\Ree(3) \cong \PGammaL(2,8)$, and $\PGammaL(2,8)$ has socle $\PSL(2,8)$. Therefore, the case $q=3$ follows from  \cite[Theorem 4.22]{BH}. Henceforth, let $q = 3^m$, where $m \geq 3$ is odd, and let $r=3^{\frac{m+1}{2}}$. Define the field automorphism $\vartheta$ of $\FF_q$ by  $x^{\vartheta}=x^{3^{\frac{m-1}{2}}}$ for every $x\in\FF_q$. Then $(x^{\vartheta})^3=x^r$ and $(x^{\vartheta})^r=x$ for every $x\in\FF_q$. For integers $i$ and $j$, we use the following mixed-exponent convention whenever the right-hand side is defined: $u^{i\vartheta+j}:=(u^{\vartheta})^iu^j.$
  Following \cite{HB,KL},  the simple Ree group   $\Ree(q)$,  of  order $q^3(q^3+1)(q-1)$,  may be realised as a subgroup of $\GL(7,q)$.
  For $x,y,z\in\FF_q$, define the following elements of   $\GL(7,q)$:
   \vspace{-5pt}$$\tau=]-1,-1,-1,-1,-1,-1,-1[, \quad  \eta=[-1,1,-1,1,-1,1,-1],$$
\vspace{-5pt}$$\alpha(x)=\small{\begin{pmatrix}
1&x^{\vartheta}&0&0&-x^{3\vartheta+1}&-x^{3\vartheta+2}&x^{4\vartheta+2}\\
0&1&x&x^{\vartheta+1}&-x^{2\vartheta+1}&0&-x^{3\vartheta+2}\\
0&0&1&x^{\vartheta}&-x^{2\vartheta}&0&x^{3\vartheta+1}\\
0&0&0&1&x^{\vartheta}&0&0\\
0&0&0&0&1&-x&x^{\vartheta+1}\\
0&0&0&0&0&1&-x^{\vartheta}\\
0&0&0&0&0&0&1\\
\end{pmatrix}},$$
$$\beta(y)=\small {\begin{pmatrix}
1&0&-y^{\vartheta}&0&-y&0&-y^{\vartheta+1}\\
0&1&0&y^{\vartheta}&0&-y^{2\vartheta}&0\\
0&0&1&0&0&0&y\\
0&0&0&1&0&y^{\vartheta}&0\\
0&0&0&0&1&0&y^{\vartheta}\\
0&0&0&0&0&1&0\\
0&0&0&0&0&0&1\\
\end{pmatrix}},
 \gamma(z)=\small {\begin{pmatrix}
1&0&0&-z^{\vartheta}&0&-z&-z^{2\vartheta}\\
0&1&0&0&-z^{\vartheta}&0&z\\
0&0&1&0&0&z^{\vartheta}&0\\
0&0&0&1&0&0&-z^{\vartheta}\\
0&0&0&0&1&0&0\\
0&0&0&0&0&1&0\\
0&0&0&0&0&0&1\\
\end{pmatrix}}.$$
   Set $\chi(x,y,z):=\alpha(x)\beta(y)\gamma(z)$  and $Q:=\{\chi(x,y,z)\mid x,y,z\in\FF_q\}$. Then $Q$ is a Sylow $3$-subgroup of  $\Ree(q)$ ($q\ge 27$), and its multiplication is given by
\vspace{-5pt}$$\chi (x_1,y_1,z_1)\chi (x_2,y_2,z_2)=\chi (x_1+x_2,y_1+y_2-x_1x_2^r,z_1+z_2-x_2y_1+x_1x_2^{r+1}-x_1^2x_2^{r}).$$
A direct calculation gives
\vspace{-5pt}$$\chi(x,y,z)^{-1}=\chi(-x,-y-x^{r+1},-z-xy-2x^{r+2}).$$
For every $k\in\FF_q^*$, define $h(k):=[k^{\vartheta},k^{1-\vartheta},k^{2\vartheta-1},1,k^{1-2\vartheta},k^{\vartheta-1},k^{-\vartheta}]$
and set $H:=\{h(k)\mid k\in \FF_q^*\}$.
In particular, $\eta=h(-1)$.
Let $a,x,y,z\in\FF_q$, let $k\in \FF_q^*$ and set $\omega:=\chi(x,y,z)$. Then the following identities hold:
\vspace{-5pt}\begin{eqnarray}\label{(0,a,0)}
 \begin{array}{lll} &&\chi(x,y,z)^{h(k)}=\alpha(xk^{r-2})\beta(yk^{1-r})\gamma(zk^{-1})=\chi(xk^{r-2},yk^{1-r},zk^{-1}),\\
 && \chi(0,a,0)^{\omega \, h(k)}=\chi(0,a,-ax)^{h(k)}=\chi(0,ak^{1-r},-axk^{-1}),\\
  &&\chi(0,0,a)^{\omega \, h(k)}=\chi(0,0,a)^{h(k)}=\chi(0,0,k^{-1}a),\\
&&(\chi(x,y,z))^3=\chi(0,0,-x^{r+2}).\end{array}
  \end{eqnarray}

\begin{prop} {\rm \cite{LN}\cite[Theorem 4.2]{Wilson} } \label{MRee}
Let $m\geq 3$ be odd, and set $q=3^m$ and $T=\Ree(q)$. Let $f$ denote the standard field automorphism of $T$ induced by $x\mapsto x^3$. Then up to conjugacy in $T$, the maximal subgroups of $T$ are as follows:
\begin{enumerate}
\item[\rm(i)]  Soluble cases: $Q:H$, where $|Q|=q^3$ and $H\cong \ZZ_{q-1}$;
 $(\ZZ_2^2\times \D_{\frac{q+1}{2}})\rtimes\ZZ_3$;
 $\ZZ_{q+r+1}\rtimes\ZZ_6$;
  $\ZZ_{q-r+1}\rtimes\ZZ_6$; and
\item[\rm(ii)] Insoluble cases:  $\ZZ_{2}\times \PSL(2,q)$; and  $\Ree(q')$, where $q'=3^l$ and $e:=\frac{m}{l}$ is prime and $q'\geq 3$.
\end{enumerate}
Here $\D_n$ denotes the dihedral group of order $n$.
\end{prop}

\begin{prop}\label{property} Let $T=\Ree(q)$, where $q=3^m$ and $m\geq 3$ is odd. With the notation introduced above, the following statements hold.
\begin{itemize}
\item[\rm(i)]Let $N_2:=\{\chi(0,y,0)\mid y\in\FF_q\}$. Then $\C_T(\eta)=\lg   \eta\rg\times \lg N_2,N_2^{\tau}\rg\cong\ZZ_2\times \PSL(2,q)$. See $\cite[sec.1; sec.2(4)]{LN}.$
\item[\rm(ii)] We have $\Z(Q)=\{\chi(0,0,z)\mid z\in \FF_q\}\cong\ZZ_3^m$. See  \cite[p.~20]{LN}.
\item[\rm(iii)] For every nontrivial subgroup $S\leq Q$, $\N_{T}(S)\leq \N_T(Q)=Q\rtimes H$.   See \cite[Lemma~1(a)]{LN}.
\item[\rm(iv)] All cyclic subgroups $L\cong \ZZ_d$ of $T$ are conjugate, for $d$ a divisor of $q-1$ with $d\geq2$;  and  $\N_{T}(L)\cong\ZZ_2\times \D_{q-1}$ and $\ZZ_2\times \PSL(2,q)$ if $d\gneqq 2$
    or $d=2$, resp.,  see \cite[sec.2(6) Lemma 2]{LN}.
\item[\rm(v)] The elements  $\chi(0,0,1)$, $\chi(0,1,0)$ and $\chi(0,-1,0)$ form a complete set of representatives for the $T$-conjugacy classes of elements of order 3. (This follows from (ii) and (iii), together with Eq(\ref{(0,a,0)}).)
\item[\rm(vi)]  Every   Sylow $2$-subgroup $S$ of $T$ is isomorphic to $\ZZ_2^3$ and
   any two $2$-subgroups having the same order are conjugate in $T$. Moreover, $\C_T(S)=S$ and $|\N_T(S)|=168$. See \cite[sec.2 (2),(3)]{LN}.
 \item[\rm(vii)] For every subgroup $L\leq T$ of order $4$ there exist a cyclic Hall subgroup $A_1$ of order $\frac{q+1}{4}$ and an element $s$ of order $6$ such that $\N_{T}(L)=\N_{T}(A_1)=L\rtimes(A_1\rtimes\lg s\rg)\cong(\ZZ_2^2\times \D_{\frac{q+1}{2}})\rtimes\ZZ_3$ and $\C_{T}(A_1)=L\times A_1$. See \cite[sec.2 (5)]{LN}.
 \item[\rm(viii)] The group $T$ has a cyclic Hall subgroup $A_0$ of order $\frac{q-1}{2}$. The group $\N_{T}(A_0)$ is dihedral of order $2(q-1)$. See \cite[sec.2 (6)]{LN}.
\item[\rm(ix)] The group $T$ has abelian Hall subgroups $A_2$ and $A_3$ of order $q-r+1$ and $q+r+1$, respectively. For each $i\in\{2,3\}$, the group $\N_{T}(A_i)$ is a Frobenius group with Frobenius kernel $A_i$ and a cyclic noninvariant factor of order $6$. See \cite[sec.2 (7)]{LN}.

\item[\rm(x)] For each $i\in\{0,1,2,3\}$ and every nontrivial subgroup $A\leq A_i$, $\N_{T}(A)\leq \N_{T}(A_i)$. Here,  $A_i$ is a cyclic subgroup of order $\frac{q-1}{2}$, $\frac{q+1}{4}$, $q-r+1$ or $q+r+1$. See \cite[sec.2 (8)]{LN}.
\end{itemize}
\end{prop}

\subsection{Actions with $M_0\cong\ZZ_2\times\PSL(2,q)$}
 As in Subsection~3.1, let $q=3^m$, where $m\geq 3$ is an odd integer and $r=3^{\frac{m+1}{2}}$. Let $T=\Ree(q)$ and  $M_0=\C_T(\eta)=\lg \eta\rg\times\lg N_2,N_2^{\tau}\rg \le T$, recalling $ \lg
 N_2,N_2^{\tau}\rg\cong \PSL(2,q)$,
\vspace{-5pt}$$\begin{array}{l}
\eta=[-1,1,-1,1,-1,1,-1],  \tau:=]-1, -1,-1,-1,-1,-1,-1[  \\
 N_2=\{ \chi(0,y,0)\mid y\in \FF_q\}\, ({\rm Proposition~\ref{property}(i)}).
 \end{array}$$
Let $f\in \Aut(T)$ be a field automorphism of order $m$. The automorphism $f$ fixes $\eta$ and hence normalizes $M_0=\C_T(\eta)$.
The following lemma gives a property about some  subgroups in $M_0$.
\begin{lem}\label{Ree-point}
Let $T=\Ree(q)$ and $M_0=\C_T(\eta)\cong\ZZ_2\times\PSL(2,q)$, where $q=3^m$ and $m\geq 3$ is odd. Choose an involution $x\in\lg N_2,N_2^\tau\rg\setminus\lg \tau\rg$ such that $x\tau=\tau x$.
Consider the action of $T$ on $\Omega_0:=[T:M_0]$. Then the following statements hold.
\begin{enumerate}
\item[\rm(1)] We have $\Fix(\lg \eta\rg)=\Fix(\lg \tau\rg)=\Fix(\lg \eta\tau\rg)=q^2-q+1$.
\item[\rm(2)] The subgroups $V_1:=\lg \eta,\tau\rg$
,$V_2:=\lg \tau,x\rg$ and $V_3:=\lg \eta\tau,\eta x\rg$ are all isomorphic to $\ZZ_2^2$ and form a complete set of representatives for the $M_0$-conjugacy classes of subgroups of $M_0$ isomorphic to $\ZZ_2^2$. Moreover, $\Fix(V_i)=q+4$ for each $i\in\{1,2,3\}$.
\item[\rm(3)] For every subgroup $S\leq M_0$ with $S\cong \ZZ_2^3$, we have $\Fix(S)=7.$
\item[\rm (4)] For every chain of subgroups $A\leq L\leq B\leq M_0$ such that $A\cong\ZZ_{\frac{q+1}{4}}$ and $B\cong\ZZ_2\times \D_{q+1}$, we have $\Fix(L)=3$.
\item[\rm(5)] For every subgroup $U\leq M_0$ with $U\cong \ZZ_3^m$, and every nontrivial subgroup $L\leq U$, we have $\Fix(L)=q$.
\item[\rm(6)] For every subgroup $L\leq M_0$ with $L\cong\ZZ_2\times \ZZ_3^m$, we have $\Fix(L)=1.$
\item[\rm(7)] For every cyclic subgroup $A_0\leq M_0$ of order $\frac{q-1}{2}$, and every subgroup $L\leq A_0$ with $|L|>2$, we have $\Fix(L)=1$.

\end{enumerate}
\end{lem}
\demo
(1) By Proposition \ref{property}(iv), the involutions of $T$ form a single $T$-conjugacy class. The subgroups of $M_0$ of order $2$ split into three $M_0$-conjugacy classes, represented by $\lg \eta\rg$, $\lg\tau\rg$ and $\lg\eta\tau\rg$. Therefore, Proposition~\ref{man} gives
\vspace{-5pt}{\small$$\begin{array}{ll}
&\Fix(\lg\eta\rg)=\Fix(\lg\tau\rg)=\Fix(\lg\eta\tau\rg)=\frac{|\N_T(\lg\eta \rg)|}{|\N_{M_0}(\lg \eta \rg)|}+\frac{|\N_T(\lg \tau \rg)|}{|\N_{M_0}(\lg \tau \rg)|}+\frac{|\N_T(\lg \eta\tau \rg)|}{|\N_{M_0}(\lg \eta\tau \rg)|}\\
&=\frac{q(q^2-1)}{q(q^2-1)}+\frac{q(q^2-1)}{2\times (q+1) }+\frac{q(q^2-1)}{2\times (q+1) }=q^2-q+1.
\end{array}$$}
~~(2) By Proposition \ref{property}(vi), all Klein four subgroups of $T$ are $T$-conjugate. We first determine the $M_0$-conjugacy classes of the Klein four subgroups contained in $M_0$.
For convenience, write $R:=\lg N_2,N_2^\tau\rg\cong\PSL(2,q),$ so that $M_0=\lg\eta\rg\times R.$
Let $K\leq M_0$ be a Klein four subgroup. We distinguish the following three cases.
First, suppose that $\eta\in K$. Then $K=\lg \eta,u\rg$ for some $u\in R$. Since all involutions of $R$ are $R$-conjugate, we may take $u=\tau$. Hence $K$ is $M_0$-conjugate to $V_1=\lg \eta,\tau\rg$.
Next, suppose that $\eta\notin K$ and $K\leq R$. Since $q\equiv 3\pmod 8$, every Klein four subgroup of $R$ is a Sylow $2$-subgroup of $R$. Hence all such subgroups are $R$-conjugate. Therefore, $K$ is $M_0$-conjugate to $V_2=\lg\tau,x\rg.$
Finally, suppose that $\eta\notin K$ and $K\nleq R$. In this case, exactly one of the three involutions of $K$ belongs to $R$, while the other two have the form $\eta u$ and $\eta v$ for two distinct commuting involutions $u,v\in R$. Their product is $\eta u\eta v=uv\in R$. Thus $\lg u,v\rg\cong \ZZ_2^2$. After conjugating in $R$, we may assume that $\lg u,v\rg=\lg \tau,x\rg$. Note that $\N_{R}(\lg \tau,x\rg)\cong \A_4$, and the elements of $\N_{R}(\lg \tau,x\rg)$ of order $3$ cyclically permute the three involutions $\tau$, $x$ and $\tau x$. It follows that all the subgroups arising in this third case are $M_0$-conjugate. We may therefore take $V_3=\lg \eta x,\eta \tau\rg$ as a representative.

Consequently, $V_1$, $V_2$ and $V_3$ form a complete set of representatives for the $M_0$-conjugacy classes of Klein four subgroups of $M_0$.
It follows from Proposition \ref{man} and Proposition \ref{property}(vii) that
{\small$
\Fix(V_i)=\frac{|\N_T(V_1)|}{|\N_{M_0}(V_1)|}
+\frac{|\N_T(V_2)|}{|\N_{M_0}(V_2)|}+\frac{|\N_T(V_3)|}{|\N_{M_0}(V_3)|}
=\frac{6(q+1)}{|\ZZ_2\times \D_{q+1}|}+\frac{6(q+1)}{|\ZZ_2\times \A_{4}|}+\frac{6(q+1)}{|\ZZ_2\times \D_{4}|}=q+4,$}
for each $i\in\{1,2,3\}.$
\vskip 3mm
(3) Let $S\leq M_0$ with $S\cong\ZZ_2^3$. Then $S$ is a Sylow $2$-subgroup of both $T$ and $M_0$. By Proposition~\ref{property}(vi), $|\N_T(S)|=168$, while $\N_{M_0}(S)\cong \ZZ_2\times \A_4$. Therefore, Proposition~\ref{man} yields
{\small$\Fix(S)=\frac{|\N_T(S)|}{|\N_{M_0}(S)|}=\frac{168}{|\ZZ_2\times \A_4|}=7.$}
\vskip 3mm
(4) Since $M_0=\C_T(\eta)$, the map $M_0g\mapsto\eta^g$ identifies the fixed points of every subgroup $H\leq M_0$ with the involutions in $\C_T(H)$.

All subgroups of $M_0$ isomorphic to $\ZZ_2\times \D_{q+1}$ are $M_0$-conjugate. Set $V:=\Z(B)\cong\ZZ_2^2$.
The subgroup $A$ is the unique Hall $2'$-subgroup of $B$.
By Proposition~\ref{property}(vii),
 $\C_T(A)=V\times A$. Then we have $\Fix(A)=3$, as $\C_T(A)$ contains exactly three involutions. Since $V=\Z(B)$, we have $\Fix(B)\geq 3$. On the other hand, $A\leq B$, and hence $\Fix(B)\leq \Fix(A)=3$. Thus $\Fix(B)=3$. This implies that $\Fix(L)=3$, as $\Fix(B)\leq\Fix(L)\leq \Fix(A)$.
\vskip 3mm
(5) Let $U\leq M_0$ with $U\cong\ZZ_3^m$. Up to $M_0$-conjugacy, we may assume that $U=N_2=\{\chi(0,t,0)\mid t\in\FF_q\}$. By Eq(\ref{(0,a,0)}) and Proposition \ref{property}(iii), we have $|\N_T(U)|=q^2(q-1)$ and $|\N_{M_0}(U)|=|\ZZ_2\times\ZZ_3^m\rtimes\ZZ_{\frac{3^m-1}{2}}|=q(q-1)$.
Therefore, Proposition~\ref{man} yields that $\Fix(U)=\frac{q^2(q-1)}{q(q-1)}=q.$
Choose $A\leq L$ with $A\cong \ZZ_3$. Up to $M_0$-conjugacy, we may assume that $A=\{\chi(0,t,0)\mid t\in\FF_3\}.$
By Proposition \ref{property}(iii) and Eq(\ref{(0,a,0)}), we get $\N_T(A)=\lg\chi(0,y,z)\mid y,z\in\FF_q\rg\rtimes \lg h(-1)\rg$ and $\N_{M_0}(A)=\lg\eta\rg\times U$. This implies that $|\N_T(A)|=2q^2$ and $|\N_{M_0}(A)|=2q$. Then by Proposition \ref{man}, we have $\Fix(A)=\frac{2q^2}{2q}=q.$
Since $A\leq L\leq U$, $\Fix(U)\leq \Fix(L)\leq \Fix(A)$. We have $\Fix(L)=q$.
\vskip 3mm
(6) Let $L\leq M_0$ with $L\cong\ZZ_2\times \ZZ_3^m,$
and let $U$ be the unique Sylow $3$-subgroup of $L$. Then $U$ is characteristic in $L$ and has order $q$. Since all Sylow $3$-subgroups of $M_0$ are conjugate, we may assume that $U=N_2$. Let $c$ be the unique involution of $L$. Since $L$ is abelian and $\C_{M_0}(U)=\lg \eta\rg\times U$, we have $\eta=c$. Therefore, $L=\lg \eta\rg\times U$.
Since $\C_T(L)\leq\C_{M_0}(U)=L$ and $L\leq \C_T(L)$, we have $\C_T(L)=L$.
 Since $\C_T(L)=L$ and $L$ contains exactly one involution, it follows that $\Fix(L)=1$.
\vskip 3mm
(7) By Proposition \ref{property}(iv), we have $\N_T(L)\cong\ZZ_2\times \D_{q-1}$. Note that $\N_{M_0}(L)\cong \ZZ_2\times \D_{q-1}$.
Then by Proposition \ref{man}, we have $\Fix(L)=\frac{2(q-1)}{2(q-1)}=1.$
\qed

\begin{lem}\label{D4}
Retain the notation of Lemma~\ref{Ree-point}, and set $B:=\lg\eta\rg\times\C_{R}(\tau)\cong\ZZ_2\times\D_{q+1}$, where $R=\lg \{\chi(0,b,0)\mid b\in\FF_q\},\tau\rg\cong\PSL(2,q).$
Let $x$ be the involution chosen in Lemma~\ref{Ree-point}, and define $E_1:=\lg \eta,\tau\rg$, $E_2:=\lg \eta,x\rg$, $E_3:=\lg \eta,\tau x\rg$, $E_4:=\lg \tau,x\rg$, $E_5:=\lg \tau,\eta x\rg$, $E_6:=\lg\eta \tau,x\rg$ and $E_7:=\lg \eta \tau,\eta x\rg.$ Then $E_1,\cdots,E_7$ form a complete set of representatives for the $B$-conjugacy classes of Klein four subgroups of $B$. Moreover, among the subgroups contained in $B$, the $M_0$-conjugacy classes represented by $V_1,V_2,V_3$ in Lemma~\ref{Ree-point} decompose as follows: $E_1,E_2,E_3\in V_1^{M_0}$, $E_4\in V_2^{M_0}$ and $E_5,E_6,E_7\in V_3^{M_0}$.
\end{lem}
\demo
Set $S:=\lg \eta,\tau,x\rg\cong\ZZ_2^3$, which is a Sylow $2$ subgroup of $B$.
Hence every Klein four subgroup of $B$ is $B$-conjugate to a subgroup of $S$, and $E_1,\cdots,E_7$ are precisely the seven Klein four subgroups of $S$.

We first identify their $M_0$-conjugacy classes. Since all involutions of $R$ are conjugate in $R$ and $\eta\in\Z(M_0)$, the subgroups $E_1,E_2,E_3$ are conjugate in $M_0$. Moreover, $\N_R(E_4)\cong A_4$ contains an element $g$ of order $3$ that cyclically permutes $\tau,x,\tau x$. Since $g$ centralizes $\eta$, it also  cyclically permutes $E_5,E_6,E_7$. As $E_1=V_1$, $E_4=V_2$ and $E_7=V_3$, this proves the displayed assignment to the three $M_0$-classes.

It remains to show that these seven subgroups are pairwise nonconjugate in $B$. Put $a:=\frac{q+1}{4}$. Then $a$ is odd, and we may write $\C_R(\tau)=\lg u,x\mid u^{2a}=x^2=1,xux=u^{-1}\rg$ and $\tau=u^a$.
Now $\Z(B)=\lg \eta,\tau\rg,$  $E_1\cap \Z(B)=E_1$, $E_2\cap \Z(B)=E_3\cap \Z(B)=\lg \eta\rg$, $E_4\cap \Z(B)=E_5\cap \Z(B)=\lg \tau\rg$ and $E_6\cap \Z(B)=E_7\cap \Z(B)=\lg \eta\tau\rg$. Since $\Z(B)$ is fixed pointwise under conjugation by $B$, only the pairs $(E_2,E_3)$, $(E_4,E_5)$ and $(E_6,E_7)$ could be $B$-conjugate.
The elements $u^ix\in\C_T(\tau)$, where $0\leq i<2a$,
form two conjugacy classes, according as $i$ is even or odd. Since $a$ is odd, $x$ and $\tau x=u^ax$ belong to different $B$-conjugacy classes. Thus $E_2$ is not $B$-conjugate to $E_3$, and $E_6$ is not $B$-conjugate to $E_7$. Finally, $E_4\leq R$, whereas $E_5\nleq R$. Since $R$ is normalized by $B$, the subgroups $E_4$ and $E_5$ are not $B$-conjugate. Therefore, the seven subgroups $E_1,\cdots,E_7$ represent distinct $B$-conjugacy classes.
\qed

\begin{lem}
Retain the notation of Lemmas~\ref{Ree-point} and~\ref{D4}. In particular, $M_0=\lg\eta\rg\times R$, where $R\cong\PSL(2,q)$.
Identify $\Omega_0=[T:M_0]$ with the $T$-conjugacy class of $\eta$, as above. Let $\alpha:=\eta$, choose $\beta\in\Omega_0\setminus \{\alpha\}$, and set $K:=T_{\alpha}\cap T_{\beta}$. If $K\neq 1$, then either $K$ contains an involution or $K=U$ for some $U\in\Syl_3(R).$
\end{lem}
\demo Suppose that $K\neq 1$ and $K$ contains no involution. Then $|K|$ is odd, and hence $K\leq R$.
Suppose that a prime $\ell\neq3$ divides $|K|$, and choose $g\in K$ of order $\ell$. Since $q\equiv 3\pmod 8$, either $\ell\mid\frac{q-1}{2}$ or $\ell\mid \frac{q+1}{4}$. In the first case, Lemma~\ref{Ree-point}(7) gives $\Fix(\lg g\rg)=1$, whereas $g\in K=T_{\alpha}\cap T_{\beta}$, a contradiction.
In the second case, choose a cyclic subgroup $A_1$ of order $\frac{q+1}{4}$ containing $\lg g\rg$. By Proposition~\ref{property}(vii) and (x), the involutions in $\C_T(g)$ are precisely the three nonidentity elements of a subgroup $V\cong\ZZ_2^2$, where $\C_T(A_1)=V\times A_1$.
Indeed, $\N_T(\lg g\rg)\leq\N_T(A_1)$, and every involution of $\N_T(A_1)\setminus\C_T(A_1)$  inverts $A_1$ and therefore does not centralize $g$. Since $g\in K$, both $\alpha$ and $\beta$ belong to $\C_T(g)$. Thus $\alpha,\beta\in V$ and $V=\lg \alpha,\beta\rg\leq\C_T(\alpha)\cap \C_T(\beta)=K$, a contradiction.
Therefore, $K$ is a nontrivial $3$-subgroup of $R$. Choose $U\in\Syl_3(R)$ with $K\leq U$. By Lemma~\ref{Ree-point}(5), $\Fix(U)=\Fix(K)$. In particular, $U$ fixes both $\alpha$ and $\beta$, so $U\leq K$. Together with $K\leq U$, this yields $K=U.$
\qed

\subsubsection{$b(G)=2$ for $\soc(G)=\Ree(q)$ and $M\cap \soc(G)\cong\ZZ_2\times\PSL(2,q)$}
The fact $b(T)=2$ was previously established by H$\acute{e}$thelyi, Horv$\acute{a}$th and Pet$\acute{e}$nyi \cite{HH}. For arbitrary almost simple extensions $T\leq G\leq \Aut(T)$, Burness, Liebeck and Shalev \cite{BL} proved the upper bound $b(G)\leq 3$. This subsection strengthens these results by proving that $b(G)=2$ for every such extension with $M\cap T\cong\ZZ_2\times \PSL(2,q).$
Set $\Gamma_0(\alpha):=\{\beta\in\Omega_0\setminus\{\alpha\}\mid T_{\alpha}\cap T_{\beta}=1\},$ where $\Omega_0=[T:M_0].$
Equivalently, $\Gamma_0(\alpha)$ is the neighborhood of $\alpha$ in the Saxl graph $\Sigma(T)$, or the union of all regular suborbits of $T$ relative to $\alpha$.
We first determine the size of $\Gamma_0(\alpha)$.
\begin{lem}
Retain the notation introduced above. In particular,
 $q=3^m\geq 27$, where $m$ is odd. Then
\vspace{-5pt}{\small$$|\Gamma_0(\alpha)|=\frac{q(q^2-1)(q+3)}{3}<\frac{|\Omega_0|}{2}.$$}
\end{lem}
\demo
Take $\alpha:=M_0\in\Omega_0$. Choose $U\in\Syl_3(R)$ and set $(K_1,K_2,K_3):=(\lg\eta\rg,\lg\tau\rg,\lg\eta\tau\rg)$,
$(K_4,K_5,K_6):=(\lg \eta,\tau\rg,\lg\tau,x\rg,\lg \eta \tau,\eta x\rg)$, and $(K_7,K_8,K_9):=(\lg\eta,\tau,x\rg,B,U)$.
By the preceding lemma, together with Lemmas~\ref{Ree-point} and \ref{D4}, every nontrivial two-point stabilizer is $M_0$-conjugate to one of these subgroups. For subgroups $A, C\leq M_0$, write $A\preccurlyeq_{M_0}C$ if an $M_0$-conjugate of $A$ is contained in $C$. The relevant subconjugacy relations are
\vspace{-5pt}\begin{equation}\label{eq1}
K_1,K_2,K_3\preccurlyeq_{M_0}K_4;\,\,
K_2\preccurlyeq_{M_0} K_5,K_6;\,\, K_3\preccurlyeq_{M_0}K_6;\,\,
K_4,K_5,K_6\preccurlyeq_{M_0}K_7\preccurlyeq_{M_0}K_8.
\end{equation}
Moreover, $K_8$ is maximal in $M_0$, $K_9\not\preccurlyeq_{M_0}K_i$ and $K_i\not\preccurlyeq_{M_0}K_9$
for any $1\leq i\leq 8$.
Together with Lemma~\ref{D4}, relations~\ref{eq1} determine which terms can occur in the fixed-point equations below.
Let $x_i$ be the number of suborbits whose point stabilizer is $M_0$-conjugate to $K_i$, and put $k_i:=\Fix_{\Omega_0}(K_i)$.
For $H,K\leq M_0$, let $f_H(K)$ denote the number of points fixed by $K$ in the action of $M_0$ on $[M_0:H]$.
It is clear that
\vspace{-5pt}\begin{equation}\label{eq2}
k_i=1+\sum\limits_{j=1}^9x_jf_{K_j}(K_i)\quad(1\leq i\leq 9).
\end{equation}
Lemma~\ref{Ree-point} gives $(k_1,\cdots,k_9)=(q^2-q+1,q^2-q+1,q^2-q+1,q+4,q+4,q+4,7,3,q).$
Using the normalizers in Proposition~\ref{property} and the conjugacy-class decompositions in Lemma~\ref{D4}, solving Eq(\ref{eq2}) successively in the order $K_9,K_8,\cdots,K_1$ gives $(x_1,\cdots,x_9)=(0,\frac{q-3}{2},\frac{q-3}{2},0,\frac{q-3}{6},\frac{q-3}{2},0,2,1)$.
Finally, we get
\vspace{-5pt}$$\begin{array}{l}
|\Gamma_0(\alpha)|=|\Omega_0|-1-\sum\limits_{i=1}^9x_i\frac{|M_0|}{|K_i|}
=\frac{q^3(q^3+1)(q-1)}{q(q^2-1)}-1-2\times\frac{q-3}{2}\times \frac{q(q^2-1)}{2}-\frac{q-3}{6}\frac{q(q^2-1)}{4}\\
-\frac{q-3}{2}\frac{q(q^2-1)}{4}-2\frac{q(q^2-1)}{2(q+1)}-\frac{q(q^2-1)}{q}
=\frac{q(q^2-1)(q+3)}{3}<\frac{|\Omega_0|}{2}
\end{array}$$
\qed

\begin{lem}\label{pslfield}
Retain the notation introduced above. Every cyclic subgroup $C\leq M$ of prime order $\ell$ with $C\nleq M_0$ is $M$-conjugate to
 $ \lg f^{\frac m{\ell}}\rg$.
\end{lem}
\demo Suppose that $\lg z\rg\cong \ZZ_{\ell}$ and $z\in M\setminus M_0$, where $M_0=\lg \eta\rg\times\lg \tau, N_2\rg \le T$.
Now, we can write $z=uf'$ with some $u\in M_0$ and $f'\in\lg f\rg$. Since $|z|=\ell$, we have that $\ell=|f'|$ is an odd prime.
It follows from $M=M_0:\lg f\rg \cong (\PSL(2,q)\times \ZZ_2):\ZZ_m$ and $\eta^f=\eta$ that $u\in \lg N_2,N_2^{\tau}\rg\cong \PSL(2,q)$
and $\lg f'\rg=\lg f^{\frac{m}{\ell}}\rg$.  By Proposition~\ref{fieldconjugate}, there exists $g\in \PGL(2,q)$ such that $(uf')^g=f'$. Since $q\equiv 3\pmod 4$, the element represented by $[-1,1]$ lies outside $\PSL(2,q)$ and commutes with $f'$. Multiplying $g$ by this element if necessary, we may assume that $g\in\PSL(2,q)$, while still having $(uf')^g=f'$. Therefore, $\lg uf'\rg$ is conjugate to $\lg f'\rg$ in $M$.
 \qed
\vskip 3mm
We now introduce  two parameters:

\vskip 3mm
{\it Let $n'(\ell)$ be the number of regular $T$-suborbits which are setwise stabilized by some subgroup $K\leq M$ of prime order $\ell$ with $K\nleq M_0$; and let $n'$ denote the total number of distinct regular $T$-suborbits setwise stabilized by some subgroup $K\leq M$ with $K\nleq M_0$.}

\vskip 3mm
 Clearly, $|\G(\alpha)|=|\G_0(\alpha)|- n'|M_0|$. Therefore, Theorem~\ref{mainpsl} holds if $|\G(\alpha)|\gneqq 0$,
\vskip 3mm

\begin{theorem}\label{mainpsl}
Let $G$ be an almost simple group with socle $T=\Ree(q)$, where $q=3^m\geq 27$ and $m$ is odd. Let $f$ be a field automorphism of $T$ of order $m$, so that $T\leq G\leq \Aut(T)=T\rtimes \lg f\rg.$ Let $M$ be a maximal subgroup of $G$ with $M\cap T=M_0\cong\ZZ_2\times \PSL(2,q)$.
 Then  the primitive action of $G$ on $\Omega=[G:M]$ by right multiplication has base size $2$.
\end{theorem}
\demo
By Remark~\ref{big}, we only need to prove that $b(G)=2$ for $G=\Ree(q):\lg f\rg$.
Set $m=r_1^{k_1} r_2^{k_2}\cdots r_{l}^{k_l}$ where $3\le r_1\lneqq r_2\lneqq \cdots \lneqq r_l$ are primes,   $m_i=\frac{m}{r_i}$, $q_i=3^{m_i}$ and $f_i=f^{m_i}.$   Then $m\ge r_1\cdots r_l\ge 3^l$ so that $l\le \log_3 m$.
By Proposition~\ref{man}, we have
\vspace{-5pt}{\small$$\Fix(\lg f_i\rg)=\frac{|\N_G(\lg f_i \rg)|}{|\N_M(\lg f_i \rg)|}=\frac{|\Ree(q_i)\rtimes\lg f\rg|}{|\ZZ_2\times\PSL(2,q_i)\rtimes\lg f\rg|}=\frac{q_i^3(q_i^3+1)(q_i-1)}{q_i(q_i^2-1)}=q_i^2(q_i^2-q_i+1).$$}In fact,  $\lg f_i\rg $ fixes exactly  $|\N_{M_0\lg f_i\rg}(\lg f_i\rg):\lg f_i\rg|=q_i(q_i^2-1)$ points on each regular $T$-suborbit  setwise stabilized by the subgroup $\lg f_i\rg$. Therefore,
\vspace{-5pt}{\small$$n'(r_i)\leq \frac{\Fix(\lg f_i\rg )}{q_i(q_i^2-1)}=   \frac{q_i(q_i^2-q_i+1)}{q_i^2-1}\leq q_i.$$}
Hence,
\vspace{-5pt}{\small$$n'\leq\sum_{i=1}^ln'(r_i)\leq\sum_{i=1}^l q_i\leq lq_1\leq3^{\frac{m}{3}}\log_3 m=q^{\frac{1}{3}}\log_3m.$$}
Therefore,
\vspace{-5pt}{\small$$
\begin{array}{l}
|\Gamma(\alpha)|=|\G_0(\alpha)|- n'|M_0|\geq\frac{q(q^2-1)(q+3)}{3}-q(q^2-1)n'\geq\frac{q(q^2-1)(q+3)}{3}-q(q^2-1)q^{\frac{1}{3}}\log_3m\\
=q(q^2-1)(\frac{q+3}{3}-q^{\frac{1}{3}}\log_3m)>0.
\end{array}$$}
This implies that $b(G)=2$.
\qed

\subsubsection{Ree unital}
A \emph{unital} is a set of $n^3 + 1$ points arranged into subsets of size $n + 1$ so that every pair of distinct points of the set is contained in exactly one subset. This is equivalent to saying that a unital is a 2-$(n^3 + 1, n + 1, 1)$ block design.

A family of unitals based on Ree groups was constructed by L\"uneburg \cite{L66}. We record, without proof, several properties of the Ree unital that will be used below (see \cite{L66}). Let $\mathcal{P}$ be the set of Sylow $3$-subgroups of $T=\Ree(q)$. Then $|\mathcal{P}|=q^3+1$ and $T$ acts $2$-transitively on $\mathcal{P}$ by conjugation.
For each involution $\mu\in T$, define $B_{\mu}:=\{P\in \mathcal{P}\mid P^{\mu}=P\},$ and set $\mathcal{B}:=\{B_{\mu}\mid \mu\in T, o(\mu)=2\}.$ Then $(\mathcal{P},\mathcal{B})$ is a $2-(q^3+1,q+1,1)$ design, called the \emph{Ree unital of order $q$}. Consequently, every two distinct points lie in a unique block, two distinct blocks meet in at most one point, and exactly $q^2$ blocks pass through each point. Moreover, $T_{B_{\mu}}=\C_T(\mu)\cong \ZZ_2\times \PSL(2,q),$ and the induced action of $T_{B_{\mu}}$ on $B_{\mu}$ is $3$-homogeneous. Since $M_0=\C_T(\eta)=T_{B_{\eta}},$  the action of $T$ on $\Omega_0$ may be identified with its action on the blocks of the Ree unital. Thus the vertices of the Saxl graph $\Sigma(T)$ may be identified with the blocks of the Ree unital, and two distinct blocks $B_{\mu_1}$ and $B_{\mu_2}$ are adjacent if and only if $T_{B_{\mu_1}}\cap T_{B_{\mu_2}}=\C_{T}(\mu_1)\cap\C_{T}(\mu_2)=1.$
\vskip 3mm
Although the Sylow $3$-subgroup of $T$ is often presented using upper unitriangular matrices, we adopt here the coordinate model used in \cite[Appendix A]{Bush} or \cite[Section 9.2.4]{Thas}. The two realizations describe equivalent seven-dimensional representations and differ only by a change of basis, a reparametrization of the root group, and the convention of left versus right action. The advantage of the model of Busch et al. is that the finite unital point is obtained directly from the fifth column of the corresponding root-group matrix. Consequently, the unital points and the standard block can be read directly from the matrix entries. This makes the subsequent block-intersection calculations more transparent. The change of model does not affect any group-theoretic conclusions.
Many of the following descriptions are adopted from \cite[Appendix A]{Bush}.

Recall that $q=3^m$. Define $\omega: \FF_q\rightarrow\FF_q,~~\omega: x\mapsto x^{3^{\frac{m+1}{2}}}.$ We identify $\mathcal{P}$ with the following set of points of the projective $6$-space $\PG(6,q)$
 represented by the following homogeneous coordinate vectors $\lg(1,0,0,0,0,0,0)\rg$ and
$$\lg(f_1(a,a',a''),-a',-a,-a'',1,f_2(a,a',a''),f_3(a,a',a''))\rg,$$
$a,a',a''\in\FF_q$ and
\vspace{-5pt}{\small$$
\left\{
  \begin{array}{ll}
    f_1(a,a',a'')=-a^{2\omega+4}-aa''^{\omega}+a^{\omega+1}a'^{\omega}+a''^2+a'^{\omega+1}-a'a^{\omega+3}-a^2a'^2, & \hbox{} \\
f_2(a,a',a'')=-a^{\omega+3}+a'^{\omega}-aa''+a^2a' ,& \hbox{} \\
f_3(a,a',a'')=-a^{2\omega+3}-a''^{\omega}+a^{\omega}a'^{\omega}+a'a''+aa'^2 .& \hbox{}
  \end{array}
\right.
$$}Set $\infty=\lg(1,0,0,0,0,0,0)\rg$, $O=\lg(0,0,0,0,1,0,0)\rg$ and $Q^*\in\Syl_3(T)$, and let $Q^*\leq T_{\infty}$. Then the Sylow $3$-subgroup $Q^*_{\infty}$ of $T$ fixing $\infty$ is given by the set of linear transformations $(a,a',a'')_{\infty}$ (acting on the left) with
$$(a,a',a'')_{\infty}:=\small{\begin{pmatrix}
1&s_1&s_2&a''&f_1(a,a',a'')&a'-a^{\omega+1}&a\\
0&1&a^{\omega}&0&-a'&0&0\\
0&0&1&0&-a&0&0\\
0&-a&a'-a^{\omega+1}&1&-a''&0&0\\
0&0&0&0&1&0&0\\
0&a^2&s_3&a&f_2(a,a',a'')&1&0\\
0&-a''-aa'&s_4&-a'&f_3(a,a',a'')&-a^{\omega}&1\\
\end{pmatrix}},$$
where
{\small$$
\left\{
  \begin{array}{ll}
    s_1=a^{\omega+3}-a'^{\omega}-aa''-a^2a', & \hbox{} \\
    s_2=-a^{2\omega+3}+a''^{\omega}+a^{\omega}a'^{\omega}+a'a''-aa'^2-a^{\omega+2}a'-a^{\omega+1}a'', & \hbox{} \\
    s_3=a^{\omega+2}+a''-aa', & \hbox{} \\
    s_4=-a^{\omega}a''+a'^2-a^{\omega+1}a'. & \hbox{}
  \end{array}
\right.
$$}Actually, all the unital points $\mathcal{P}$ can be viewed as column-vectors $\infty:=(1,0,0,0,0,0,0)^T$
and the orbit of $Q^*_{\infty}$ on $O=(0,0,0,0,1,0,0)^T$,
that is, $p(a,a',a''):=(a,a',a'')_{\infty}O$ and so
\vspace{-5pt}\begin{eqnarray}\label{unitalpoint}
\begin{array}{l}
p(a,a',a'')=\lg(f_1(a,a',a''),-a',-a,-a'',1,f_2(a,a',a''),f_3(a,a',a''))^T\rg.
\end{array}
\end{eqnarray}
And we deduce
\vspace{-5pt}\begin{eqnarray}\label{a}
\begin{array}{l}
(a,a',a'')_{\infty}\cdot p(b,b',b'')=p(a+b,a'+b'+a^{\omega}b,a''+b''-ab'+a'b-a^{\omega+1}b), \end{array}
\end{eqnarray}
for any $b,b',b''\in\FF_q.$
The two-point stabilizer $T_{(\infty,O)}$ consists of the elements $\varrho(t)$, acting by
\vspace{-5pt}\begin{eqnarray}\label{t}
\begin{array}{l}
\varrho(t)\cdot p(b,b',b'')=p(bt,b't^{\omega+1},b''t^{\omega+2}),~~b,b',b''\in\FF_q~~\text{and}~~t\in\FF_{q}^*;
\end{array}
\end{eqnarray}
see \cite[p. 90]{Hao}. This corrects the exponent in the first coordinate of the corresponding formula in \cite[Appendix A]{Bush}, where $bt^{\omega}$ should read $bt.$
The structure of $T$ yields $\N_{T}(Q^*)=Q^*\rtimes H^*=T_{\infty},$ where $H^*=\{\varrho(t)\mid t\in\FF_q^*\}.$
We find that $\varrho(-1)$ is an involution fixing the point set $\{\infty\}\cup\{p(0,b',0)\mid b'\in\FF_q\},$
which is--by definition--the block $B(\infty,O)$, where $B(\infty,O)$ denotes the block containing $\{\infty, O\}$. The other blocks through $\infty$ are the images of $B(\infty,O)$ under the elements $(a,0,a'')_{\infty}$, $a,a''\in\FF_q$ (Since $(0,a',0)_{\infty}$ stabilises $B(\infty,O)$). This implies that the blocks through $\infty$ are
\vspace{-5pt}\begin{eqnarray}\label{block}
\begin{array}{l}
B_{a,a''}:=\{\infty\}\cup\{p(a,b',a''-ab')\mid b'\in\FF_q\},~~a,a''\in\FF_q.
\end{array}
\end{eqnarray}
In particular, we set $B_0:=B_{0,0}.$
Set $U=\{(0,b,0)_{\infty}\mid b\in\FF_q\}$. By Eq(\ref{a}) and Eq(\ref{t}), we have
$T_{(\infty,B_0)}=U\rtimes H^*\lesssim \ZZ_2\times \PSL(2,q).$
From now on we set $A=T_{(\infty,B_0)}\cong\ZZ_3^m\rtimes\ZZ_{3^m-1}$ for convenience.
Considering the action of $A$ on the set $\{B_{a,c}\mid a,c\in\FF_q\}$, by Eq(\ref{a}), Eq(\ref{t}) and the definition of $B_{a,c}$, we have
\begin{eqnarray}\label{AB(a,c)}
\begin{array}{ll}
(0,b,0)_{\infty}\cdot B_{a,c}=B_{a,c-ab},&\varrho(t)\cdot B_{a,c}=B_{at,ct^{\omega+2}}.
\end{array}
\end{eqnarray}

\begin{lem}\label{AN}
Retain the notation introduced above.
 The three $A$-orbits on the blocks through $\infty$ are $\{B_{0,0}\}$, $\{B_{0,c}\mid c\in\FF_q^{*}\}$
  and $\{B_{a,c}\mid a\in\FF_q^*,c\in\FF_q\}$. Their block stabilizers in $A$ are $A$, $U$ and $1$, respectively.
Consequently, for every block $C\neq B_0$ through $\infty$,
\vspace{-5pt}$$\begin{array}{l}
T_{B_0}\cap T_C=\left\{
  \begin{array}{ll}
    U, & \hbox{$C=B_{0,c}$ with $c\neq 0$,} \\
    1, & \hbox{$C=B_{a,c}$ with $a\neq 0$.}
  \end{array}
\right.
\end{array}$$
\end{lem}
\demo
The block $B_{0,0}$ is fixed by $A$, by the definition of $A$.
Now take $B_{0,c}$ with $c\neq 0$. Eq(\ref{AB(a,c)}) shows that $U\leq T_{B_{0,c}}$.
Eq(\ref{AB(a,c)}) together with the fact $\{t^{\omega+2}\mid t\in\FF_q^*\}=\FF_{q}^*$, shows that $H^*$ is regular on $\{B_{0,c}\mid c\in\FF_q^*\}.$ Hence, the second set is one orbit of length $q-1$, with stabilizer $U$.
Finally, suppose $a\neq 0$. Eq(\ref{AB(a,c)}) shows that
$H^*$ is regular on the set $\{B_{a,0}\mid a\in\FF_q^*\}$.
Once $a_0$ is fixed, Eq(\ref{AB(a,c)}) shows that $U$ is regular on the set $\{B_{a_0,c}\mid c\in\FF_q\}.$
Therefore, $\{B_{a,c}\mid a\in\FF_q^*,c\in\FF_q\}$ forms one $A$-orbit. Its length is $|A|=q(q-1)$, so its block stabilizer is trivial.

Finally, for any block $C\neq B_{0}$ through $\infty$, we have $T_{B_0}\cap T_{C}\leq T_{(\infty,B_0)}$. This implies that $T_{B_0}\cap T_{C}=A_C$ for each block $C\in\mathcal{B}\setminus \{B_0\}$ with $B_0\cap C=\{\infty\}$.
The conclusion follows from the arguments above.
\qed

\begin{defi}\label{def}
Let $P\in\mathcal{P}$ be a unital point and let $B\in\mathcal{B}$ be a unital block containing $P$. Define the local bad set $\mathcal{C}_P(B):=\{C\in\mathcal{B}\mid T_B\cap T_C\neq 1, P\in B\cap C\}$ for the incident pair $(P,B).$
\end{defi}
The natural action of $T$ on the unital points is $2$-transitive, so it is transitive on the incident point-block pairs. This implies that $|\mathcal{C}_P(B)|=|\mathcal{C}_{\infty}(B_0)|=|\{B_0\}\cup\{B_{0,c}\mid c\in\FF_q^*\}|$.
Thus, by Lemma \ref{AN} $|\mathcal{C}_P(B)|=q$ for any incident pair $(P,B)$.

\subsubsection{Proof of the BG-Conjecture for the Saxl graph $\Sigma(T)$ when $T=\Ree(q)$ and $M_0\cong\ZZ_2\times\PSL(2,q)$}
Recall that the vertices of the Saxl graph $\Sigma(T)$ with $M_0\cong\ZZ_2\times\PSL(2,q)$  may be identified with the blocks of the Ree unital. Retain the notation above. We let $\Gamma_0(B)$ denote the neighborhood of $B$ in  $\Sigma(T)$.
The main result of this subsection is the following theorem.
\begin{theorem}
In the Saxl graph $\Sigma(T)$ with $M_0\cong\ZZ_2\times\PSL(2,q)$, for any two distinct unital blocks $B$ and $D$, we have
\vspace{-5pt}$$
|\Gamma_0(B)\cap \Gamma_0(D)|\geq\left\{
  \begin{array}{ll}
    q(q-2), & \hbox{$B\cap D\neq \emptyset$;} \\
    (q+1)(q-7), & \hbox{$B\cap D=\emptyset$.}
  \end{array}
\right.
$$
\end{theorem}
This theorem follows from Lemmas~\ref{BDintersect} and~\ref{BDemptyset}.

\begin{lem}\label{BDintersect}
Let $B,D\in\mathcal{B}$ be distinct blocks such that $B\cap D=\{P\}$, for some unital point $P\in\mathcal{P}.$
Then $|\Gamma_0(B)\cap \Gamma_0(D)|\geq q(q-2)$.
\end{lem}
\demo
By the preceding calculation, we have $|\mathcal{C}_P(B)\cup \mathcal{C}_P(D)|\leq 2q$. Then $|\Gamma_0(B)\cap \Gamma_0(D)|\geq q^2-2q$, as there are exactly $q^2$ blocks through a fixed unital point $P$.
\qed
\vskip 3mm
By Lemma \ref{AN}, the local bad set of the pair $(\infty, B_0)$ is $\mathcal{C}_{\infty}(B_0)=\{B_{0,c}\mid c\in\FF_q\}.$
Delete the common point $\infty$ and take the union of these $q$ blocks. Eq(\ref{block}) gives the set
\vspace{-5pt}\begin{eqnarray}\label{Pi}
\begin{array}{l}
\Pi_0:=\{p(0,s,c)\mid s,c\in\FF_q\}.
\end{array}
\end{eqnarray}
Clearly, $\{B_{0,c}\setminus\{\infty\}\mid c\in\FF_q\}$ is a  partition of $\Pi_0$.
We now consider two disjoint blocks. As preparation, we first give the following lemma.
\begin{lem}\label{PiD}
Retain the notation introduced above. Let $D$ be a unital block satisfying $D\cap B_0=\emptyset$. Then
$|\Pi_0\cap D|\leq 4.$
\end{lem}
\demo
Recall that $B_0=\{\infty\}\cup\{p(0,b,0)\mid b\in\FF_q\}$.
By Eq(\ref{unitalpoint}), we get $p(0,b,0)=\lg(b^{\omega+1},-b,0,0,1,b^{\omega},0)^T\rg$.
The map  $\phi_1:\PG(1,q)\rightarrow B_0$, defined by
$$\phi_1: \left\{
                                         \begin{array}{ll}
                                           \lg(x,y)\rg=\lg(\frac{x}{y},1)\rg\mapsto p(0,\frac{x}{y},0)=\lg((\frac{x}{y})^{\omega+1},-\frac{x}{y},0,0,1,(\frac{x}{y})^{\omega},0)^T\rg, & \hbox{$y\neq 0$;} \\
                                           \lg(x,y)\rg=\lg(1,0)\rg\mapsto \infty, & \hbox{y=0,}
                                         \end{array}
                                       \right.
$$ is a bijection.
Now, we may identify $B_0$ with the projective line $\PG(1,q),$
 and we rewrite $B_0=\{\lg(1,0)\rg\}\cup\{\lg(b,1)\rg\mid b\in\FF_q\}.$
Let $D=gB_0$ be a block such that $B_0\cap D=\emptyset$,
for some $g\in T$. Consequently, the point of $D$
corresponding to $\lg (x,y)\rg$ is $Q_{x,y}:=\lg gv(x,y)\rg$, where $v(x,y):=(x^{\omega+1},-xy^{\omega},0,0,y^{\omega+1},x^{\omega}y,0)^T.$
Let $\mathcal{H}_0:=\{\lg(x_1,x_2,x_3,\cdots,x_7)\rg\in\PG(6,q)\mid x_3=0\}.$
Since the third coordinate of $p(a,a',a'')$ is $-a$ (by Eq(\ref{unitalpoint})), the definition of $\Pi_0$ gives
 $\mathcal{H}_0\cap \mathcal{P}=\{\infty\}\cup \Pi_0.$
Then for any point $Q_{x,y}\in D=gB_0$, we  have
\vspace{-5pt}$$Q_{x,y}\in D\cap \Pi_0\Leftrightarrow Q_{x,y}\in\mathcal{H}_0\Leftrightarrow (gv(x,y))_3=0,$$
where $(x_1,x_2,x_3,\cdots,x_7)_3$ denotes the third coordinate of the vector $(x_1,x_2,x_3,\cdots,x_7)$ for some $x_i\in\FF_q$ with $1\leq i\leq 7.$
This implies that $|D\cap \Pi_0|=|\{\lg (x,y)\rg\in\PG(1,q)\mid (gv(x,y))_3=0\}|.$
Write $g=(g_{ij})_{1\leq i,j\leq 7}\in T$. Then
\begin{eqnarray}\label{eq}
\begin{array}{l}
(gv(x,y))_3=g_{31}x^{\omega+1}-g_{32}xy^{\omega}+ g_{35}y^{\omega+1}+g_{36}x^{\omega}y.
\end{array}
\end{eqnarray}
The solutions of
Eq(\ref{eq}) are the projective points $\lg (x,y)\rg\in\PG(1,q).$
This implies that
$|D\cap \Pi_0|=|\{\lg(x,y)\rg\in\PG(1,q)\mid g_{31}x^{\omega+1}-g_{32}xy^{\omega}+ g_{35}y^{\omega+1}+g_{36}x^{\omega}y=0\}|$.
\vskip 3mm
Now, we will prove $|D\cap \Pi_0|\leq 4$.
Set $F(x,y)=g_{31}x^{\omega+1}-g_{32}xy^{\omega}+ g_{35}y^{\omega+1}+g_{36}x^{\omega}y$.
Then it is sufficient to get the number of solutions $\lg (x,y)\rg\in\PG(1,q)$ of $F(x,y)=0$ for $x,y\in\FF_q$.
If $F(x,y)=0$ has at most two projective solutions, there is nothing to prove.
Now, suppose that it has at least three distinct projective solutions.
Note that $T_{D}\cong\ZZ_2\times\PSL(2,q)$ and is $3$-homogeneous on the unital points of $D$.
So, without loss of generality, we assume that $\lg(1,0)\rg$, $\lg(0,1)\rg$ and $\lg(1,1)\rg$ are three projective solutions of $F(x,y)=0$. This implies that
$g_{31}=0$, $g_{35}=0$ and $g_{31}-g_{32}+g_{35}+ g_{36}=0.$
Finally, we have $F(x,y)=g_{32}(x^{\omega}y-xy^{\omega}).$

\vskip 3mm
{\bf Case 1: $g_{32}=0$}
\vskip 3mm
In this case, all the points $\lg(x,y)\rg\in\PG(1,q)$ are solutions. This implies that $D\subset \Pi_0$. Note that $\{B_{0,c}\setminus\{\infty\}\mid c\in\FF_q\}$ partitions $\Pi_0$ into $q$ disjoint parts. Since $q+1$ unital points of $D$ all lie in $\Pi_0$, this forces  two distinct points of $D$ to lie in the same block $B_{0,c_0}$ for some $c_0\in\FF_q$.
Since Ree unital is a $2-(q^3+1,q+1,1)$ design, we have $D=B_{0,c_0},$ with $B_{0,c_0}=\{\infty\}\cup\{p(0,b,c_0)\mid b\in\FF_q\}$. So, $\infty\in D\cap B_0$, a contradiction to $D\cap B_0=\emptyset$.
\vskip 3mm
{\bf Case 2: $g_{32}\neq0$}
\vskip 3mm
In this case $|D\cap \Pi_0|=|\{\lg(x,y)\rg\in\PG(1,q)\mid F(x,y)=0\}|=|\{\lg(x,y)\rg\in\PG(1,q)\mid x^{\omega}y=xy^{\omega}\}|.$
The case $y=0$ gives the solution $\lg (1,0)\rg$. If $y\neq 0$, setting $z:=\frac{x}{y}$ gives $z^{\omega}=z$.
Recall that $\omega: \FF_q\rightarrow\FF_q,~~\omega: x\mapsto x^{3^{\frac{m+1}{2}}}.$ Since $(\frac{m+1}{2},m)=1$, we have that the fixed field of $\omega$ is the set $\FF_3$. Thus, the projective solutions of $F(x,y)=0$ are $\lg (1,0)\rg$, $\lg (0,1)\rg$
$\lg (1,1)\rg$ and $\lg (-1,1)\rg$. This completes the proof.
\qed

\begin{lem}\label{BDemptyset}
Let $B,D$ be two unital blocks satisfying $B\cap D=\emptyset$.
Then, in the Saxl graph $\Sigma(T)$, we have $|\Gamma_0(B)\cap \Gamma_0(D)|\geq (q+1)(q-7)$.
\end{lem}
\demo
For each pair $(P_B,P_D)\in B\times D$, let $C_{P_B,P_D}$ be the unique unital block containing $P_B$ and $P_D$, and set $\mathcal{X}:=\{C_{P_B,P_D}\mid P_B\in B, P_D\in D\}.$ We first show that these blocks are pairwise distinct. Suppose that $C_{P_B,P_D}= C_{P_B',P_{D}'}$.
If $P_B\neq P_B'$, then the block $C_{P_B,P_D}$ contains two distinct points $P_B$ and $P_B'$. This implies that $C_{P_B,P_D}=B$ and $P_D\in B\cap D$, a contradiction. Hence, $P_B=P_B'$. By the same argument, $P_D=P_D'$. Therefore, $|\mathcal{X}|=(q+1)^2$.
Fix $P_B\in B$. Since $T$ is transitive on the incident point-block pairs, choose $g\in T$ such that $g\cdot(P_B,B)=(\infty, B_0)$. Then $(g\cdot D)\cap B_0=\emptyset$ and $|\Pi_0\cap (g\cdot D)|\leq 4$ by Lemma~\ref{PiD}.
By Lemma~\ref{AN}, the block $C_{P_B,P_D}$ is not adjacent to $B$ if and only if $g\cdot P_D\in \Pi_0$.
Thus, for each fixed $P_B\in B$, at most four choices of $P_D\in D$ give blocks $C_{P_B,P_D}$ that are not adjacent to $B$ in the Saxl graph $\Sigma(T)$.
Consequently, at most $4(q+1)$ members of $\mathcal{X}$ are not adjacent to $B$. By symmetry, at most $4(q+1)$ members of $\mathcal{X}$ are not adjacent to $D$. Therefore, $|\Gamma_0(B)\cap \Gamma_0(D)|\geq (q+1)^2-8(q+1)=(q+1)(q-7).$
\qed

\subsection{Actions with $M_0\cong\Ree(q')$}
Retain the notation introduced above. Recall that $T = \Ree(q)$, where $q=3^m\geq 27$, and let $f$ be a field automorphism of $T$ of order $m$. Let $G$ be an almost simple group satisfying $T\leq G\leq \Aut(T)$. Write $m = m'e$, where $e\geq3$ is an odd prime, and set $q'=3^{m'}$. Let $M$ be a maximal subgroup of $G$ such that $M_0=M\cap T=\Ree(q')$.
Recall that $r=3^{\frac{m+1}{2}}$ and set $r'=3^{\frac{m'+1}{2}}$. Then \vspace{-5pt}$$q^2-q+1=(q+r+1)(q-r+1),\,
 q'^2-q'+1=(q'+r'+1)(q'-r'+1).$$
 By Remark \ref{big}, it suffices to take $G=T\rtimes\lg f\rg$ and $M=M_0\rtimes\lg f\rg$. In what follows,  we shall determine the representatives $\lg z\rg$  for the $M$-conjugacy classes  of prime-order subgroups contained in $M_0$ in Lemma~\ref{M0} and those not contained in  $M_0$ in Lemma~\ref{M}, respectively.
 For each representative $\lg z\rg$, we define its contribution to $\check{Q}(G)$ by
\vspace{-5pt}$$ {\small Q(\lg z \rg):=\frac{|M|}{|\Omega|}\frac{(|\lg z\rg|-1)\Fix(\lg z \rg)}{|\N_{M}(\lg z\rg)|}.}$$
 Using the estimates obtained in Lemmas~\ref{M0} and~\ref{M}, we prove that $\check{Q}(G)<\frac{1}{2}$ in Lemma~\ref{mainfield} and then the BG-Conjecture in this case is proved by Proposition~\ref{prob}.

\begin{lem}\label{M0}
Retain the notation introduced above. Every subgroup of $M_0$ of prime order $\ell$ is  $M$-conjugate to $\lg z\rg $,  with the following possibilities:
 \begin{enumerate}
\item[\rm(i)]$\ell=2$: $\lg\eta\rg$,   with $Q(\lg \eta \rg)\leq0.0175$;
\item[\rm(ii)]  $\ell=3$: $\lg\chi(0,0,1)\rg $ or $ \lg\chi(0,1,0) \rg$,   with  $Q(\lg  \chi(0,0,1) \rg)+Q(\lg  \chi(0,1,0) \rg)<0.0185$;
\item[\rm(iii)] $\ell\mid \frac{q'+1}4$: a single $M$-conjugacy class represented by $\lg z\rg $ for each fixed prime $\ell$, with  $Q(\lg  z\rg)\leq\frac{98}{2025q'^{6e-13}}$;
\item[\rm(iv)] $\ell\mid (q'\pm r'+1)$: a single $M$-conjugacy class represented by $\lg z\rg $ for each fixed prime $\ell$,  with $Q(\lg  z\rg)\leq\frac{14q'^{14}(q+r+1)}{75q^7(q'\pm r'+1)}$;
\item[\rm(v)]  $\ell\mid \frac{q'-1}2$:a single $M$-conjugacy class represented by $\lg z\rg $ for each fixed prime $\ell$, with  $Q(\lg  z\rg)\leq\frac{21}{50q'^{6e-13}}$.
\end{enumerate}
\end{lem}
\demo
First, since $q'\geq3$ and $q\geq27$, we compute
$$\begin{array}{l}
\frac{|M|}{|\Omega|}
=\frac{mq'^6(q'^3+1)^2(q'-1)^2}{q^3(q^3+1)(q-1)}\leq\frac{28mq'^{14}}{25q^7}.
\end{array}$$
By Propositions ~\ref{MRee} and~\ref{property}, we collect the required data for the prime-order subgroups $\lg z\rg$ of $M_0$ in the following cases. In each case, the subgroups of $M_0$ that are $G$-conjugate to the chosen subgroup form a single $M_0$-conjugacy class, so Proposition~\ref{man} gives the fixed-point formulas below.

(i) By Proposition~\ref{property}(iv), $M_0$ has a single conjugacy class of involutions, so we may take $z=\eta$.
Since   $\N_G(\lg \eta \rg)\cong (\ZZ_2\times\PSL(2,q))\rtimes\lg f\rg$ and $\N_M(\lg \eta \rg)\cong (\ZZ_2\times\PSL(2,q'))\rtimes\lg f\rg$, we have $\Fix(\lg \eta \rg)=\frac{|(\ZZ_2\times \PSL(2,q))\rtimes\lg
f\rg |} {|(\ZZ_2\times\PSL(2,q'))\rtimes\lg f\rg |}=\frac{q(q^2-1)}{q'(q'^2-1)}$  and then
{\small$$\begin{array}{l}
Q(\lg \eta \rg)=\frac{|M|}{|\O|}\frac{(|\lg \eta \rg|-1)\Fix(\lg \eta \rg)}{|\N_M(\lg \eta \rg)|}=\frac{|M|}{|\O|}\frac{q(q^2-1)}{q'^2(q'^2-1)^2m}\leq\frac{28mq'^{14}}{25q^7}\frac{q^3}{q'^6(1-q'^{-2})^2m}\leq\frac{567}{400q'^{4e-8}}\leq0.0175.\end{array}$$}

(ii) By Proposition~\ref{property}(v), there are two $M_0$-conjugacy classes of subgroups of order $3$, represented by
 $\lg \chi(0,0,1) \rg $ and $\lg \chi(0,1,0) \rg$, respectively.

Let $z=\chi(0,0,1)$. It follows from Eq(\ref{(0,a,0)}) that
\vspace{-5pt}$$\begin{array}{lll}\N_G(\lg z \rg)&=&\lg\chi (a,b,c),h(k),f\mid a,b,c\in\FF_q, k=\pm 1 \rg\cong (Q\rtimes\ZZ_{2})\rtimes\lg f \rg ,\\
\N_M(\lg z \rg)&=&\lg\chi (a,b,c),h(k),f\mid a,b,c\in\FF_{q'}, k=\pm 1\rg\cong (Q_{q'}\rtimes\ZZ_{2})\rtimes\lg f \rg , \end{array} $$
where $Q_{q'}\in \Syl_3(\Ree(q'))$.
Then $\Fix(\lg z\rg )=\frac{q^3}{q'^3}$ and {\small$Q(\lg z\rg)=\frac{|M|}{|\O|}\frac{q^3}{q'^6m}<0.01383$.}
Let $z=\chi(0,1,0)$. It follows from Eq(\ref{(0,a,0)}) that
\vspace{-5pt}$$\begin{array}{lll}\N_G(\lg z \rg)&=&\lg\chi (0,b,c),h(k)\mid b,c, k\in\FF_q, k^{2-2r}=1\rg \rtimes\lg f \rg ,\\
\N_M(\lg z \rg)&=&\lg\chi (0,b,c),h(k)\mid b,c, k\in\FF_{q'}, k^{2-2r'}=1\rg\rtimes\lg f \rg .\end{array} $$
Then  $\Fix(\lg z\rg)=\frac{q^2}{q'^2}$ and {\small$Q(\lg z\rg)=\frac{|M|}{|\O|}\frac{q^2}{q'^4m}<0.00461$.}

Consequently, we have $Q(\lg  \chi(0,0,1) \rg)+Q(\lg  \chi(0,1,0) \rg)<0.0185$.
\vskip 3mm

(iii) Suppose  $\ell\mid \frac{q'+1}{4}$.
 Since $\Ree(q')$ has a single conjugacy class of cyclic subgroups of order $\frac{q'+1}{4}$, we may choose  $\lg z\rg $ in any such cyclic subgroup.
  By Propositions \ref{MRee} and \ref{property}(vii), we have $\N_G(\lg z\rg)\cong \ZZ_2^2\times \D_{\frac{q+1}{2}}\rtimes\ZZ_3\rtimes\lg f \rg$ and $\N_M(\lg z \rg)\cong \ZZ_2^2\times \D_{\frac{q'+1}{2}}\rtimes\ZZ_3\rtimes\lg f \rg$.
Hence $\Fix(\lg  z\rg)=\frac{q+1}{q'+1}$ and
{\small$Q(\lg z\rg)=\frac{|M|}{|\O|}\frac{(\ell-1)(q+1)}{6(q'+1)^2m}\leq\frac{|M|}{|\O|}\frac{(q+1)}{24q'm}\leq\frac{98}{2025q'^{6e-13}}.$}

\vskip 3mm
(iv)
Suppose  $\ell\mid q'\pm r'+1$.
 Since $\Ree(q')$ has a single conjugacy class of cyclic subgroups of order  $q'\pm r'+1$, we may choose  $\lg z\rg $ in any such cyclic subgroup.
By Propositions \ref{MRee} and \ref{property}(ix,x), we have $|\N_G(\lg z\rg)|\leq 6(q+r+1)m$ and $|\N_M(\lg z \rg)|= 6(q'\pm r'+1)m$. Hence,
\vspace{-5pt}$$\begin{array}{lll}
Q(\lg z\rg)&\leq& \frac{|M|}{|\O|}\frac{(\ell-1)\frac{q+r+1}{q'\pm r'+1}}{6(q'\pm r'+1)m}=\frac{|M|}{|\O|}\frac{(\ell-1)(q+r+1)}{6(q'\pm r'+1)^2m}\leq\frac{|M|}{|\O|}\frac{(q+r+1)}{6(q'\pm r'+1)m}\leq\frac{14q'^{14}(q+r+1)}{75q^7(q'\pm r'+1)}.
\end{array}$$
\vskip 3mm
(v) Suppose  $\ell\mid \frac{q'-1}{2}$.  Since $\Ree(q')$ has a single conjugacy class of cyclic subgroups of order  $\frac{q'-1}{2}$, we may choose  $\lg z\rg $ in any such cyclic subgroup.
Then $|\N_G(\lg  z\rg)|=2m(q-1)$ and $|\N_M(\lg  z\rg)|=2m(q'-1)$.
Hence $\Fix(\lg  z\rg)=\frac{q-1}{q'-1}$ and
{\small$
Q(\lg z\rg)=\frac{|M|}{|\O|}\frac{(\ell-1)(q-1)}{2(q'-1)^2m}\leq\frac{21}{50q'^{6e-13}}.
$}
\qed

\begin{lem} \label{M}
Retain the notation introduced above. Every subgroup of $M$ of prime order $\ell$ that is not contained in $M_0$ is  $M$-conjugate to $\lg z\rg $,  with the following possibilities:
 \begin{enumerate}
\item[\rm(i)]  $\ell\neq e$: $\lg f'\rg$, where $f'=f^{\frac m\ell}$, with $Q(\lg  f'\rg)\leq\frac{63(\ell-1)}{25q'^{14/3}}$;
\item[\rm(ii)] $\ell=e$, and either $q'\geq 27$, or $q'=3$ and $m=e\geq5$: $M$-conjugacy classes with representatives $\lg a\rg$ satisfying $ \sum_{\lg a \rg} Q(\lg a\rg)\leq \frac{56e m^2q'^{19}}{75q^7};$
\item[\rm(iii)] $\ell=e$ with $q'=3$ and $m=e=3$: $M$-conjugacy classes with representatives $\lg a\rg$ satisfying $ \sum_{\lg a \rg} Q(\lg a\rg)\leq 0.16.$
\end{enumerate}
In (ii) and (iii), each sum runs over one representative from each $M$-conjugacy class of subgroups of order $e$ in $M$ that are not contained in $M_0$.
\end{lem}
\demo
Recall from the preceding lemma that $\frac{|M|}{|\O|}\leq \frac{28mq'^{14}}{25q^7}$.
Let $z\in M\setminus M_0$ have prime order $\ell$. Replacing $z$ by a suitable generator of $\lg z\rg$, we may write $z=uf'$, where $u\in M_0$ and $f'=f^{\frac{m}{\ell}}$. In what follows, we consider the cases $\ell=e$ and $\ell\neq e$ separately.

\vskip 3mm
(1)
Suppose that $\ell\neq e$. Then $\ell\mid m'$. Write $m=m'e=m_1\ell e$.  The restriction of $f'=f^{m_1e}$ to $M_0$ is a field automorphism of order $\ell$, with fixed-point subgroup $\Ree(3^{m_1})$.
By Proposition \ref{fieldconjugate}, $z=f^{m_1e}$ up to conjugacy in $M_0$.
Hence, $\N_G(\lg z \rg)\cong \Ree(3^{\frac{m}{\ell}})\rtimes\lg f \rg$ and $\N_M(\lg z \rg)\cong \Ree(3^{\frac{m'}{\ell}})\rtimes\lg f \rg$.
Note that $|\Ree(q)|=q^7(1+q^{-3})(1-q^{-1})$.
Then we have $\frac{2}{3}q^7\leq|\Ree(q)|\leq q^7$.
Therefore,
\vspace{-5pt}$$\begin{array}{lll}
Q(\lg z\rg)=\frac{|M|(\ell-1)}{|\O|}\frac{|\Ree(q^{1/\ell})|}{m|\Ree(q'^{1/\ell})|^2}
\leq\frac{|M|(\ell-1)}{|\O|}\frac{q^{7/\ell}}{m(\frac{2}{3}q'^{7/\ell})^2}\leq\frac{63(\ell-1)}{25}\frac{q'^{14-14/\ell}}{q^{7-7/\ell}}
\leq\frac{63(\ell-1)}{25q'^{14/3}}.
\end{array}$$

(2)  Suppose that $\ell=e$.  In this case, $f'=f^{m'}$ centralizes $M$, so $(uf')^\ell=1$ if and only if $u^\ell=1$. Then the $M$-conjugacy class of the subgroup $\lg uf'\rg$ is determined by that of $u$.
By Proposition~\ref{fieldconjugate}, all subgroups of order $\ell$ in $M$ not contained in $M_0$ are $G$-conjugate to $\lg f'\rg$. Moreover, $\N_G(\lg f'\rg)=\N_M(\lg f'\rg)=M$.
In the following sums, $K$ runs over one representative from each $M$-conjugacy class of these subgroups. By Proposition~\ref{man}, we have
$\Fix(\lg f'\rg)=\Fix(\lg uf'\rg)=\sum\limits_{K}\frac{|\N_G(K)|}{|\N_M(K)|}=|M|\sum\limits_{K}\frac{1}{|\N_M(K)|}$.
Consequently,
\vspace{-5pt}$$
\begin{array}{l}
\sum\limits_{K}Q(K)=\frac{|M|}{|\Omega|}(\ell-1)\Fix(\lg f'\rg)\sum\limits_K\frac{1}{|\N_M(K)|}=\frac{(\ell-1)(\Fix(\lg f'\rg))^2}{|\Omega|}.
\end{array}
$$
Thus it suffices to bound $\Fix(\lg f'\rg)$.
If $\ell\nmid |M_0|$, then $u=1$ and $\Fix(\lg f'\rg)=1$.
It remains to consider the following four cases: $\ell=3$, $\ell\mid \frac{q'+1}{4}$, $\ell\mid q'\pm r'+1$ and $\ell\mid \frac{q'-1}{2}$.
\vskip 2mm
(2.1)  $\ell=3$
\vskip 2mm
 Since $M_0$ contains exactly  three classes of elements of order $3$ with representatives $ u_i $ ($i=1, 2, 3$), where  $ u_1=\chi(0,0,1),$ $ u_2=\chi(0,1,0) $ and  $u_3=\chi(0,-1,0)$, the group $M$  contains exactly four $M$-conjugacy classes of subgroups of order $3$ not contained in $M_0$, represented by $\lg u_if'\rg $ ($i=1, 2, 3$) and $\lg f'\rg$.
Since $\N_G(\lg f'\rg)=\N_M(\lg f'\rg)=\Ree(q'):\lg f\rg$,  $\N_M(\lg u_1f'\rg)\geq\lg\chi(x,y,z),h(1),f\mid x,y,z\in\FF_{q'}\rg$,
$\N_M(\lg
u_2 f'\rg)\geq\lg\chi(0,y,z),h(1),f\mid y,z\in\FF_{q'}\rg$
and $\N_M(\lg u_3 f'\rg)\geq\lg\chi(0,y,z),h(1),f\mid y,z\in\FF_{q'}\rg$, we get
\vspace{-5pt}$$\begin{array}{lll}
 \Fix(\lg f'\rg)&=&\Fix(\lg u_if'\rg)\leq
\frac{q'^3(q'^3+1)(q'-1)m}{q'^3m}+\frac{q'^3(q'^3+1)(q'-1)m}{q'^2m}+\frac{q'^3(q'^3+1)(q'-1)m}{q'^2m}\\&+&\frac{q'^3(q'^3+1)(q'-1)m}{q'^3(q'^3+1)(q'-1)m}=(q'^3+1)(q'-1)(2q'+1)+1<2.42q'^5+1. \end{array}$$
~~~~(2.2)  $\ell\mid\frac{q'+1}{4}:$
\vskip 2mm
It follows from Proposition \ref{property}(vii) that $\N_{M_0}(\lg u \rg)\cong\ZZ_2^2\times \D_{\frac{q'+1}{2}}\rtimes\ZZ_3$ and $\C_{M_0}(u)\cong\ZZ_2^2\times \ZZ_{\frac{q'+1}{4}}$. This implies that there are $\frac{\ell-1}{6}$ $M_0$-conjugacy classes of elements of order $\ell$. Every subgroup of $M$ of order $\ell$ not contained in $M_0$ is conjugate to either $\langle f' \rangle$ or $\langle u f' \rangle$ for some $u \in M_0$ of order $\ell$. Therefore, the number of $M$-conjugacy classes of subgroups of order $\ell$ not contained in $M_0$ is at most $1 + \frac{\ell-1}{6}$. Furthermore, $\langle f' \rangle$ and $\langle u f' \rangle$ are conjugate in $G$, as follows from Proposition \ref{fieldconjugate}. Hence, $\Fix(\lg
f'\rg)=\Fix(\lg u f'\rg)\leq \frac{|\Ree(q'):\lg f\rg|}{|\Ree(q'):\lg f\rg|}+\frac{\ell-1}{6}\frac{|\Ree(q'):\lg f\rg|}{(q'+1)e}\leq 1+\frac{q'^3(q'^2-q'+1)(q'-1)m}{6}\leq1+\frac{1}{6}mq'^6$, as $\N_M(\lg
uf'\rg)\geq \C_M(\lg uf'\rg)\gtrsim \ZZ_2^2\times \ZZ_{\frac{q'+1}{4}}\rtimes\lg f' \rg$.
\vskip 3mm
(2.3)  $\ell\mid q'\pm r'+1$
\vskip 2mm

Let $\lg u\rg$ be a subgroup of $M_0$ of order $\ell$. In this case, we have $\N_{M_0}(\lg u\rg)\cong(\ZZ_{q'\pm r'+1}\rtimes\ZZ_6)$ and $\C_{M_0}(u)\cong\ZZ_{q'\pm r'+1}$, by Proposition \ref{property}(ix). By the same argument as in case (2.2), there are at most $1+\frac{\ell-1}{6}$ $M$-conjugacy classes of subgroups of order $\ell$ that are not contained in $M_0$. Then $\Fix(\lg f'\rg)=\Fix(\lg u f'\rg)\leq
\frac{|\Ree(q'):\lg f\rg|}{|\Ree(q'):\lg f\rg|}+\frac{\ell-1}{6}\frac{|\Ree(q'):\lg f\rg|}{(q'\pm r'+1)\ell}\leq1+\frac{q'^3(q'^2-1)(q'+r'+1)m}{6}\leq1+\frac{7mq'^6}{18}$, as $\N_M(\lg uf'\rg)\geq \C_M(\lg uf'\rg)\gtrsim
\ZZ_{q'\pm r'+1} \rtimes\lg f'\rg$.
\vskip 3mm
(2.4)  $\ell\mid \frac{q'-1}{2}$
\vskip 2mm
 In this case, we have $\N_{M_0}(\lg u\rg)\cong(\ZZ_{2}\times \D_{q'-1})$ and $\C_{M_0}(u)\cong \ZZ_2\times\ZZ_{\frac{q'-1}{2}}$,  by Proposition \ref{property}(viii). By the same argument as in case (2.2), there are at most $1+\frac{\ell-1}{2}$ $M$-conjugacy classes of subgroups of order $\ell$ that are not contained in $M_0$. Then $\Fix(\lg f'\rg)=\Fix(\lg uf'\rg)\leq \frac{|\Ree(q'):\lg f\rg|}{|\Ree(q'):\lg f\rg|}+\frac{\ell-1}{2}\frac{|\Ree(q'):\lg f\rg|}{(q'-1)\ell}\leq1+\frac{q'^3(q'^3+1)m}{2}\leq1+\frac{14}{27}mq'^6$, because $\N_M(\lg uf'\rg)\geq \C_M(\lg uf'\rg)\gtrsim
\ZZ_2\times\ZZ_{\frac{q'-1}{2}}\rtimes \lg f'\rg$.

Since $m,q'\geq 3$, the preceding bounds imply, in all cases, that $\Fix(\lg f'\rg)\leq\frac{2}{3}mq'^6$.
Also, the initial estimate and $|M_0|\geq \frac{2q'^7}{3}$ give $\frac{1}{|\Omega|}=\frac{|M|}{|\Omega|}\frac{1}{m|M_0|}\leq\frac{42q'^7}{25q^7}$.
It follows that $\sum_KQ(K)=\frac{(e-1)(\Fix(\lg f'\rg))^2}{|\Omega|}\leq \frac{56e m^2q'^{19}}{75q^7}$.
Finally,
 if $q'=3$ and $m=e=3$, then the more precise estimate above gives $\Fix(\lg f'\rg)\leq 393$. Hence  $\sum_{K}Q(K)\leq \frac{2\times393^2\times 42}{25\times 3^{14}}\leq0.16$.

 This completes the proof.
\qed

\begin{lem}\label{mainfield} Let $G=\Ree(q):\lg f\rg $, where $q=3^m\geq27$ and $m$ is an odd integer. Let $M=\Ree(q'):\lg f\rg $, where $q'=3^{m'}$ and $e=\frac{m}{m'}\geq3$ is an odd prime. Consider the primitive action of $G$ on $\O:=[G:M]$ by right multiplication. Then $b(G)=2$ and the BG-Conjecture holds.
\end{lem}
\demo
Retain the notation introduced above. If $m'>1$, then
write $m'=r_1^{k_1} r_2^{k_2}\cdots r_{l}^{k_l}$, where $3\le r_1\lneqq r_2\lneqq \cdots \lneqq r_l$ are distinct primes and $k_i\geq 1$ for all $i$. Set $m_i=\frac{m}{r_i}$ and $f_i=f^{m_i}.$ Let $\lg z_{q'-1}\rg$, $\lg z_{q'+1}\rg$, $\lg z_{q'+r'+1}\rg$ and $\lg z_{q'-r'+1}\rg$ be cyclic Hall subgroups of $M_0$ of orders $\frac{q'-1}{2}$, $\frac{q'+1}{4}$, $q'+r'+1$ and $q'-r'+1$, respectively.
Further, let
\vspace{-5pt}\begin{enumerate}
  \item  [{\rm(a)}] $q'-1=2t_1^{j_1} t_2^{j_2}\cdots t_{o}^{j_o}$, where $5\le t_1\lneqq t_2\lneqq \cdots \lneqq t_o$ are distinct primes;
    \item  [{\rm(b)}] $q'+1=2^{2} s_1^{i_1} s_2^{i_2}\cdots s_{h}^{i_h}$,
 where $5\le s_1\lneqq s_2\lneqq \cdots \lneqq s_h$ are distinct primes;
  \item [{\rm(c)}] $q'+r'+1=u_1^{\ell_1} u_2^{\ell_2}\cdots u_{s}^{\ell_{s}}$, where $5\le u_1\lneqq u_2\lneqq \cdots \lneqq u_s$ are distinct primes; and
 \item [{\rm(d)}] $q'-r'+1=w_1^{d_1} w_2^{d_2}\cdots w_{t}^{d_{t}}$, where $5\le w_1\lneqq w_2\lneqq \cdots \lneqq w_t$ are distinct primes.
\end{enumerate}
Furthermore, denote by $z_{q'-1, i}$ an element of order $t_i$ in $\lg z_{q'-1}\rg$; denote by $z_{q'+1, i}$ an element of order $s_i$ in $\lg z_{q'+1}\rg$; denote
 by $z_{q'+r'+1, i}$ an element of order $u_i$ in $\lg z_{q'+r'+1}\rg$; denote
 by $z_{q'-r'+1, i}$ an element of order $w_i$ in $\lg z_{q'-r'+1}\rg$.
Recall that $\mathcal{P}^*(M):=\{\lg g_1\rg,\lg g_2\rg,\cdots,\lg g_k\rg\}$ is a set of representatives for the $M$-conjugacy classes of prime-order subgroups of $M$.

If $q'\geq27$, then we have $o\leq\log_5\frac{q'-1}{2}$, $h\leq\log_5\frac{q'+1}{4}$, $s+t\leq 2\log_5(q'+r'+1)$ and $\sum_{i=1}^l(r_i-1)\leq m'-1\leq q'.$
By Lemmas \ref{M0}, \ref{M} and Proposition \ref{prob}, we have
\vspace{-5pt}$${\small\begin{array}{lll}
 \check{Q}(G)&=&\frac{|M|}{|\O|}(\sum_{i=1}^k(|\lg g_i\rg|-1)\frac{\Fix(\lg g_i\rg)}{|\N_M(\lg g_i\rg )|})=\sum_{i=1}^kQ(\lg g_i\rg)\\
 &\leq&Q(\lg \eta \rg)+Q(\lg \chi(0,1,0)\rg)
+Q(\lg \chi(0,0,1) \rg)+\sum_{i=1}^oQ(\lg z_{q'-1,i}\rg)
+\sum_{i=1}^hQ(\lg z_{q'+1,i}\rg)\\
&&+\sum_{i=1}^sQ(\lg z_{q'+r'+1,i}\rg)+\sum_{i=1}^tQ(\lg z_{q'-r'+1,i}\rg)+\sum_{|f_i|\neq e}Q(\lg f_i\rg)+\sum_{\lg a\rg}Q(\lg a\rg)\\
&\leq& 0.036+\log_5\frac{q'-1}{2}\times \frac{21}{50q'^{6e-13}}+\log_5\frac{q'+1}{4}\times \frac{98}{2025q'^{6e-13}}+2\log_5(q'+r'+1)\times\frac{14q'^{14}(q+r+1)}{75q^7(q'-r'+1)}\\
&&+\sum\limits_{i=1}^l\frac{63(r_i-1)}{25q'^{14/3}}+\frac{56e m^2q'^{19}}{75q^7}<\frac{1}{2}.
\end{array}}$$

If $q'=3$, we have $o=h=t=0$ and $s=1$.
By Lemmas \ref{M0}, \ref{M} and Proposition \ref{prob}, we have
\vspace{-5pt}$${\small\begin{array}{lll}
 \check{Q}(G)&=&\frac{|M|}{|\O|}(\sum_{i=1}^k(|\lg g_i\rg|-1)\frac{\Fix(\lg g_i\rg)}{|\N_M(\lg g_i\rg )|})=\sum_{i=1}^kQ(\lg g_i\rg)\\
 &\leq&Q(\lg \eta \rg)+Q(\lg \chi(0,1,0)\rg)
+Q(\lg \chi(0,0,1) \rg)+Q(\lg z_{q'+r'+1,1}\rg)+\sum_{\lg a\rg}Q(\lg a\rg)\\
&\leq& 0.036+\frac{14q'^{14}}{75q^7}\frac{(q+r+1)}{(q'+ r'+1)}+0.16<\frac{1}{2}\\
\end{array}}$$
as desired.$\hskip 13.8cm \Box$

\section{Groups with $\soc(G)=\Sz(q)$}
In this section, we consider the case $\soc(G)=\Sz(q)$. The notation introduced below replaces that used specifically for the Ree case and will be used throughout the remainder of the paper.
Let $T:=\soc(G)=\Sz(q)$, where $q=2^m$ for  some  odd integer $m$. Then $T\le G\le \Aut(T)=T\rtimes\lg f \rg$, where $f$ is the field automorphism of order $m\geq 3$.
 To prove Theorem \ref{main}, we need to
consider all primitive permutation representations of $G$ on the set $\O:=[G:M]$ of right cosets of $M$ in $G$, for  maximal subgroups $M$ of $G$. The case when $M$ is soluble was settled by Burness and Huang \cite{BH} and  we
may therefore assume that $M$ is insoluble. By the maximal subgroup classification recalled in Subsection 4.1, it remains to consider the case in which $M_0:=M\cap T=\Sz(q')$, where $q'=2^l$ $(l\geq 3)$ for  some divisor $l$ of $m$ with $\frac{m}{l}$ a prime.  Thus we may write
\begin{center}
 { \it $G=T\rtimes\lg f_1\rg,$ where $|f_1|=m_1$,  for some $m_1\di m$;  $M_0=\Sz(q')$ and  $M=M_0.\lg f_1\rg .$ }\end{center}
In Subsection 4.1, some preliminary results for $\Sz(q)$ are collected;   in Subsection 4.2, the action of $T$ on $\O_0:=[T:M_0]$ is considered and the size of the union $\G_0$ of all regular suborbits of $T$ is determined; and  finally in Subsection 4.3,  Theorem~\ref{main} for $\soc(G)=\Sz(q)$ is proved.

\subsection{Preliminaries for $\Sz(q)$}
The Suzuki groups $\Sz(q)$ are an infinite family of simple groups of Lie type, first defined by Suzuki in \cite{suzuki1,suzuki2} as subgroups of $\SL(4, q)$. As our analysis necessitates explicit computations, we recall the essential definitions and theorems from \cite{suzuki2}, adopting the framework established in \cite{Fang-Praeger2}. A shorter, alternative introduction to these groups may be found in \cite{Huppert3,Wilson}.

 Let $q=2^m$, where $m\geq 3$ is an odd integer.  Set
$r:=2^{\frac{m+1}{2}}$ so that $r^2=2q$. Define $\vartheta\in \Aut(\FF_q)$ by $\a^\vartheta= \a^r$ for all $\a\in \FF_q$. Then $(\a^\vartheta)^\vartheta=\a^2$.
For $\alpha,\beta\in
\FF_q$, define
 \vspace{-5pt}$$\chi(\alpha,\beta):=\begin{pmatrix}
1&0&0&0\\
\alpha&1&0&0\\
\alpha_{31}&\alpha^r&1&0\\
\alpha_{41}&\beta&\alpha&1\\
\end{pmatrix},$$
where $\alpha_{31}:=\alpha^{1+r}+\beta$ and $\alpha_{41}:=\alpha^{2+r}+\alpha\beta+\beta^r$.
For any $k\in \FF_q^*$,  let $\kappa(k):=[k^{\frac{1}{2}r+1}, k^{\frac{1}{2}r}, k^{-\frac{1}{2}r}, k^{-\frac{1}{2}r-1}]\in \SL(4,q)$.  Let
$\tau:=]1,1,1,1[.$ Set
\vspace{-5pt}$$\begin{array}{lll} &Q:=Q(q)=\{\chi(\alpha,\beta)\mid \alpha,\beta\in\FF_q\},  &K:=K(q)=\{\kappa(k)\mid k\in\FF_q^*\}\cong \ZZ_{q-1},\\
 &H:=H(q):=\lg Q,K \rg, \, &T:=\lg H,\tau \rg.\end{array}$$
 \f Then $T$ is  isomorphic to the Suzuki simple group $\Sz(q)$, which is of order $q^2(q^2+1)(q-1)$. Note that $q^2+1=(q+ r+1)(q-r+1)$ and $5\mid q^2+1$. Given that $q + r + 1$ and $q - r + 1$ are coprime, then $5$ must divide exactly one of these two factors. Let $\delta\in\{1,-1\}$ be the unique sign such that $5 \mid q + \delta r + 1$. For each subfield $\FF_{q_1} \subseteq \FF_{q}$, define $Q(q_1)$ and $K(q_1)$ by restricting the parameters in the above matrices to $\FF_{q_1}$, with the parameter $k$ nonzero, and set $H(q_1):=\lg Q(q_1), K(q_1)\rg$ and $\Sz(q_1)=\lg H(q_1),\tau\rg$.

\begin{prop}
{\rm\cite[Theorem 9]{suzuki2}}\label{MSz}
Up to conjugacy, every  maximal subgroup of $T=\Sz(q)$ is one of the following:
\begin{enumerate}
\item[\rm(i)]  Soluble cases:  $Q.K$;  $\D_{2(q-1)}$;  $\ZZ_{q+\delta r+1}\rtimes\ZZ_4$ with $5\mid q+\delta r+1$;  $\ZZ_{q-\delta r+1}\rtimes\ZZ_4$ with $5\nmid q-\delta r+1$; and
\item[\rm(ii)] Insoluble cases: $\Sz(q')$, where $q'=2^l>2$, $l\mid m$ and $\frac{m}{l}$ is prime.
\end{enumerate}
\end{prop}

Throughout this section,  by  $A_1$ and  $A_2$ we denote cyclic subgroups of orders $q+\delta r+1$ and $q-\delta r+1$ respectively. Let $A_0=K$.
Set $B_i=\N_T(A_i)$, where $i=0, 1, 2$. Since $T$ has a single conjugacy class of subgroups isomorphic to $A_i$, we may choose $A_i$ and $B_i$ to be normalized by the field automorphism $f$. For a group $G$ we write $G^*:=G\setminus\{1_G\}$. Next, we give some properties of the subgroups of $T$.

\begin{prop} {\rm\cite{Fang-Praeger2,Huppert3,suzuki2}} \label{Sz}
 Retain the notation introduced above, and let $T=\Sz(q)$. Then the following statements hold.
\begin{description}
\item [\rm(i)]$\Aut(T)=T\rtimes\lg f\rg$;
\item[\rm(ii)] We have $Q\in \Syl_2(G)$ and $\Z(Q)=\{\chi(0,\beta)\mid \beta\in \FF_q\}\cong \ZZ_2^m$. Every element $x\in Q\setminus Z(Q)$ has
    order $4$ and is not conjugate to $x^{-1}$ in $T$;
\item[\rm(iii)] For any $l\di m$,  $Q\cap \Sz(2^l)$ is the unique Sylow $2$-subgroup of $\Sz(2^l)$ containing $\chi(0,1)$. If $l>1$, its center is $\{  \chi(0,\alpha)\di \alpha\in \FF_{2^l}\}$;
\item[\rm(iv)] We have $\N_{T}(Q)=H=Q.K$,  $H\cap H^{\tau}=K$ and $\kappa(k)^{\tau}=\kappa(k)^{-1}$ for every $k\in\FF_q^*$;
    \item[\rm(v)] The group $K$ acts regularly by conjugation on both $Z(Q)^*$ and $(Q/Z(Q))^*$;
\item[\rm(vi)] For every $x\in T$,  either  $Q\cap Q^x=1$ or $Q= Q^x$;
\item[\rm(vii)] We have $B_0=\lg K,\tau \rg \cong D_{2(q-1)}$;
\item[\rm(viii)] For $i=1,2$, $|B_i|=4|A_i|=4(q+(-1)^{i+1}\delta r+1)$ and $B_i$ is a Frobenius group;
\item[\rm(ix)] For each $i\in\{1,2\}$, we may choose an element $t_i\in B_i$ of order $4$ such that $B_i=\N_T(A_i)= A_i\rtimes \lg t_i\rg\cong\ZZ_{q\pm \d r+1}\rtimes \ZZ_4$ and $a^{t_i}=a^q$ for every $a\in A_i.$
\end{description}
\end{prop}

\f The following formulas can be verified directly; some of them also appear in \cite{suzuki2}.
\begin{prop} \label{fomula} Retain the notation introduced above. For all $\a,\b,\g,\e,x,y\in\FF_q$ and $k\in\FF_q^*$, the following identities hold.
 \begin{description}
\item[\rm(1)] $\chi(\alpha,\beta)\chi(\gamma,\e)=\chi(\alpha+\gamma,\alpha\gamma^r+\beta+\e),$ $\chi(\alpha,\beta)^{-1}=\chi(\alpha,-\alpha^{1+r}-\beta)$;
\item[\rm(2)] $\chi(\alpha,\beta)^{\chi(x,y)}=\chi(\alpha ,\alpha x^r+x\alpha^r+\beta)$;
\item[\rm(3)] $\chi(\alpha,\beta)^{\chi(x,y)\kappa(k)}=\chi(\alpha k,(\alpha x^r+x\alpha^r+\beta)k^{1+\vartheta})$;
\item[\rm(4)] $\lg \chi(\alpha,\beta) \rg=\{\chi(\alpha,\beta),\chi(0,\alpha^{1+r}),\chi(\alpha,\beta+\alpha^{1+r}), \chi(0,0)\}$;
\item[\rm(5)] $\kappa(k)\chi(\alpha,\beta)=\chi(\alpha k^{-1},\beta k^{-1-\vartheta})\kappa(k)$,
$\kappa(k)^{-1}\chi(\alpha,\beta)\kappa(k)=\chi(\alpha k,\beta k^{1+\vartheta})$.
\end{description}
\end{prop}
\subsection{The action of $T$ on $[T:M_0]$}
Let $T=\Sz(q)$ and $M_0=\Sz(q')$, where $q=2^m$ and $q'=2^l>2$, $m>1$ is odd, and $e:=m/l$ is prime. In this subsection, we shall consider the primitive action of $T$ on $\O_0=[T :M_0]$. Take $\a:=M_0\in \O_0$, so that $T_{\a}=M_0$.
The following key  lemma concerning the  structure of subgroups of $T$ is used  to determine the possible two-point stabilizers
$(M_0)_{\b}\leq M_0=T_{\a}\in \O_0$ with $\b\in \O_0 \setminus \{\a\}$.   We define $r' := 2^{\frac{l+1}{2}}$ and let $\delta' \in \{1, -1\}$ be determined by the condition $5 \mid q' + \delta' r' + 1$.

\begin{lem}\label{Sz-point}
Let $T=\Sz(q)$ and $M_0=\Sz(q')$, where $q=2^m$, $q'=2^l$, $m,l > 1$ are odd integers, and $\frac{m}{l}$ is prime. Set $Q_0 = Q(q')$ and $K_0 = K(q')$, as defined in the preceding subsection.
For the action of $T$ on $\O_0$, the following statements hold.
\begin{enumerate}
\item[\rm(1)] For any $\ZZ_4\lesssim L\leq Q_0$, we have $\Fix(L)=2^{m-l}$;
\item[\rm(2)] For any $Q_0\lneqq L\leq  Q_0.K_0$, we have
   $\Fix(L)=1$;
\item[\rm(3)] For any $\ZZ_{l'}\cong L\le M_0$, where $l'\mid q'-1$ and $l'\ne 1$, we have $\Fix(L)=\frac{q-1}{q'-1}$;
\item[\rm(4)] For any $\ZZ_{l'}\rtimes\ZZ_2\cong L\le M_0$, where $l'\mid q'-1$ and $l'\neq 1$, we have $\Fix(L)=1$;
\item[\rm(5)] For any  $\ZZ_{l'}\lesssim L\leq \ZZ_{q'\pm \delta'r'+1}\le M_0$, where  $l'\mid q'\pm \delta'r'+1$ and $l'\neq 1$, we have
 $\Fix(L)=\frac{q\pm \delta r+1}{q'\pm \delta'r'+1}$;
\item[\rm(6)] For any  $ \ZZ_{l'}\rtimes\ZZ_2\lesssim L \lesssim \ZZ_{q'\pm \delta'r'+1}\rtimes\ZZ_4\le M_0$, where $l'\mid q'\pm \delta'r'+1$ and $l'\neq1$, we have
     $\Fix(L)=1$;
\item[\rm(7)] For any  $\ZZ_2\cong L\le M_0$, we have $\Fix(L)=2^{2(m-l)}$;
\item[\rm(8)] For   $L=\Z(Q_0)\cong \ZZ_2^l$,  we have  $\Fix(L)=2^{2(m-l)}$;
\item[\rm(9)] For any $\ZZ_2^s\rtimes \ZZ_{s'}\lesssim L\le Z(Q_0):K_0$, where $s$ and $s'$ are integers such that $2\leq s\leq l$, $s'>1$ and $s'\mid 2^s-1$, we have
  $\Fix(L)=1$.
\end{enumerate}
\end{lem}
\demo
(1) Let $\ZZ_4\cong A\leq Q_0\leq Q$, where $A=\lg \chi(\alpha,\beta)\rg $ with $\alpha\neq 0$. Then $A=\{\chi(\alpha,\beta),
\chi(0,\alpha^{1+r}),\chi(\alpha,\beta+\alpha^{1+r}),\chi(0,0)\}.$

First, we show that $\N_T( A )\cong\Z(Q).\ZZ_2$. In fact,
for any $x\in \N_T(A)$,  we have $Q\cap Q^x\ne 1$  and so $Q^x=Q$ by   Proposition~\ref{Sz}.(vi). Therefore,   $\N_T(A)\leq \N_T(Q)=H$.
 For any element $u\in  \N_T(A)$, write $u=\chi(x,y)\kappa(k).$ Since $\chi(0,\alpha^{1+r})$ is the unique involution of $A$, we have $\chi(0,\alpha^{1+r})^u=\chi(0,\alpha^{1+r}).$
  By Propositions~\ref{Sz}(ii) and~\ref{fomula}(2-4),  from
 \vspace{-5pt}$$\chi(0,\alpha^{1+r})^{\chi(x,y)\kappa(k)}
=\chi(0,\alpha^{1+r}k^{1+\vartheta})=\chi(0,\alpha^{1+r}),$$
we get $k^{1+\vartheta}=k^{1+r}=1$, that is $k=1$;   and from
\vspace{-5pt}$$\chi(\alpha,\beta)^{\chi(x,y)\kappa(1)}
=\chi(\alpha ,(\alpha x^r+x\alpha^r+\beta))
\in\{\chi(\alpha,\beta),\chi(\alpha,\beta+\alpha^{1+r})\},$$
we get that $y\in \FF_q$ and either
$\alpha x^r+x\alpha^r+\beta=\beta$ or $\alpha x^r+x\alpha^r+\beta=\beta+\alpha^{1+r}$. Both  equations  give  $x\in \{ 0, \alpha\}$, meaning
 $\N_T(A)=\{\chi(0,y),\chi(-\alpha,y)\mid y\in\FF_q\}\cong\Z(Q).\ZZ_2$, as desired.
 Also, we have $\N_{M_0}(A)\cong\Z(Q(q')).\ZZ_2$.
  It follows from Proposition \ref{Sz}(v) that the cyclic subgroups of order $4$ form a single $M_0$-conjugacy class.
Hence, by Proposition~\ref{man}, we get
$\Fix(A)=|\N_T(A):\N_{M_0}(A)|=2^{m-l}$.

Since $\Z(Q)\leq\N_T(Q_0)$, we have $\Z(Q)\N_{M_0}(Q_0)\leq \N_T(Q_0)$, which implies that $\frac{|\Z(Q)|}{|\Z(Q)\cap \N_{M_0}(Q_0)|}\leq \frac{|\N_{T}(Q_0)|}{|\N_{M_0}(Q_0)|}$.
Noting that $\mathrm{Fix}(Q_0) \leq \mathrm{Fix}(A)$ and using Proposition \ref{man}, we obtain
\vspace{-5pt}$$\begin{array}{lll} 2^{m-l}&=&\frac{|Z(Q)|}{|Z(Q)\cap Q_0K_0|}\leq\Fix(Q_0)\le |\Fix(A)|=2^{m-l},\end{array}$$
 forcing equality at each step. In particular, $|\Fix(L)|=2^{m-l}$, for any subgroup $A\le L\le Q_0$. Moreover, the proof above shows that $\N_T(Q_0) = \Z(Q) \N_{M_0}(Q_0)=\Z(Q)(Q_0K_0)$.
\vskip 3mm

(2) For any $L$ with $Q_0 \lneqq L \le Q_0K_0$, we have $Q_0K_0 \le \N_{M_0}(L)$, and since $Q_0K_0$ is maximal in $M_0$, it follows that $\N_{M_0}(L) = Q_0K_0$. We shall show that $\N_T(L) = Q_0K_0$ as well; hence, by Proposition \ref{man}, $\Fix(L) = 1$, as any $T$-conjugate of $L$ contained in $M_0$ is $M_0$-conjugate to $L$.

Since $Q_0$ is a characteristic subgroup of $L$, we have $\N_T(L) \le \N_T(Q_0) = \Z(Q)\N_{M_0}(Q_0) = \Z(Q)(Q_0K_0)$. Take arbitrary elements $u \in L\setminus Q_0$ and $v \in \N_T(L)$, and write them as $u = \chi(s, t)\kappa(k_1)$ and $v = \chi(\alpha, \beta)\kappa(k_2)$ for some $k_1, k_2 \in \mathbb{F}_{q'}^*$.  Using the formula
\vspace{-5pt}$${\small\begin{array}{lll}
&&(\chi(s,t)\kappa(k_1))^{\chi(\alpha,\beta)\kappa(k_2)}\\
&=&\chi([\alpha(1+k_1^{-1})+s]k_2,[\alpha^{1+r}k_1^{-r}+s\alpha^rk_1^{-r}+\alpha s^r-\alpha^{1+r}-\beta+t+\beta k_1^{-1-\vartheta}]k_2^{1+\vartheta})\kappa(k_1),
\end{array}}$$
where $s,t,\alpha,k_1,k_2\in\FF_{q'}$ and $\b\in\FF_q$, we note that the right hand side belongs to $L$ only if all parameters are in $\FF_{q'}$. Thus $\beta \in \FF_{q'}$, and therefore $v \in Q_0K_0$. Thus we have shown $\N_T(L) =\N_{M_0}(L)= Q_0K_0$, as desired.
\vskip 3mm

(3)-(4)  Let $L\cong \ZZ_{l'}$ be a subgroup of $M_0$ with $l'\mid q'-1$ and $l'\ne 1$. Then $\N_T(L)\cong\D_{2(q-1)}$ and $\N_{M_0}(L)\cong\D_{2(q'-1)}$. Since all subgroups of $M_0$ isomorphic to $\ZZ_{l'}$ are conjugate, we obtain $\Fix(L)=|\N_T(L):\N_{M_0}(L)|=\frac{q-1}{q'-1}$, by Proposition \ref{man}. Similarly, let $L\cong \ZZ_{l'}\rtimes\ZZ_2$ be a subgroup of $M_0$, where $l'\mid q'-1$ and $l'\ne 1$. In this case, $\Fix(L)=|\N_T(L):\N_{M_0}(L)|=|L:L|=1.$

\vskip 3mm

(5)-(6)  Let $\ZZ_{l'}\cong L\leq \ZZ_{q'\pm\delta'r'+1}$. The cyclic subgroups of order $l'$ form a single $M_0$-conjugacy class.  So $\Fix(L)=|\N_T(L):\N_{M_0}(L)|=|(\ZZ_{q\pm \delta
r+1}\rtimes\ZZ_4):(\ZZ_{q'\pm \delta' r'+1}\rtimes\ZZ_4)|=\frac{q\pm \delta r+1}{q'\pm \delta' r'+1}$.
Set $\ZZ_{l'}\rtimes\ZZ_2\cong J\le A\leq \ZZ_{q'\pm \delta' r'+1}\rtimes\ZZ_4$, for $l'\mid q'\pm \delta' r'+1$. It follows from Proposition \ref{Sz}(ix) that
 $\N_{M_0}(J)\cong\ZZ_{l'}\rtimes\ZZ_4$ and $\N_T(J)\cong\ZZ_{l'}\rtimes\ZZ_4$. Note that $1\leq \Fix(A)\leq\Fix(J)$.
 Since the subgroups of $M_0$ isomorphic to $J$ form a single conjugacy class, Proposition~\ref{man} gives $\Fix(J)=|\N_T(J):\N_{M_0}(J)|=1$.
Hence,
$\Fix(A)=1$.

\vskip 3mm

(7)  Let $L\cong\ZZ_2$. Since  the involutions of $M_0$ form a single conjugacy class, it follows from Proposition~
\ref{Sz}.(ii,v,vi) that   $\Fix(L)=|\N_T(L):\N_{M_0}(L)|=|Q:Q_0|=2^{2(m-l)}$.
\vskip 3mm

(8)  Let $\ZZ_2^l\cong L=\Z(Q_0)$. The elementary abelian subgroups of order $q'$ in $M_0$ form a single conjugacy class. Hence, $\Fix(L)=|\N_T(L):\N_{M_0}(L)|=|Q.K_0:Q_0.K_0|=2^{2(m-l)}$, which
follows from Proposition~
\ref{Sz}(v,vi).

\vskip 3mm

(9) We first show that  $\N_T(L)=\N_{M_0}(L)$, so that  $\Fix(L)=1$, as the $T$-conjugates of $L$ contained in $M_0$ form a single $M_0$-conjugacy class.
 Let $E=L\cap \Z(Q_0)$, the unique Sylow $2$-subgroup of $L$. Since $E$ is characteristic in $L$, we have $\N_T(L)\leq \N_T(E)\leq \N_T(Q)$ by Proposition~\ref{Sz}(vi).
Take any element $u_1v_1 \in \N_T(L)$, and write $u_1=\chi(\alpha,\beta)$ and $v_1=\kappa(k_1)$. For any element $u_2v_2 \in L \le \Z(Q_0):K_0$, we can write $u_2=\chi(0,t)$ and $v_2=\kappa(k_2)$.
From
\vspace{-5pt}$$\begin{array}{lll} &&(\chi(0,t)\kappa(k_2))^{\chi(\alpha,\beta)\kappa(k_1)}\\
&=&\chi(\alpha(1+k_2^{-1})k_1,[\alpha^{1+r}k_2^{-r}-\alpha^{1+r}-\beta+t+\beta k_2^{-1-\vartheta}]k_1^{1+\vartheta})\kappa(k_2),\end{array}$$
we get $\alpha=0$, $k_1\in\FF_{q'}$ and $\beta\in\FF_{q'}$,  which implies  $\N_T(L)\leq\N_{M_0}(L)$. Hence, $\N_T(L)=\N_{M_0}(L)$, as desired.
\qed

\vskip 3mm

By Proposition~\ref{MSz} and  Lemma~\ref{Sz-point},  the point stabilizer of each nontrivial nonregular suborbit is $M_0$-conjugate to one of the following subgroups: $K_1=Q_0$, $K_2\cong\ZZ_{q'-1}$,
$K_3\cong \ZZ_{q'+\delta'r'+1}$, $K_4\cong \ZZ_{q'-\delta'r'+1}$ and $K_5=Z(Q_0)\cong \ZZ_2^l$, where $5\mid q'+\delta' r'+1$ and $\d'=\pm 1$. For $1\leq i\leq 5$,  let $x_i$ denote the number of suborbits whose point stabilizer is $M_0$-conjugate to $K_i$.

\begin{lem}\label{nonregular-Su}
With the notation introduced above and $l \geq 3$, we have:
\vspace{-5pt}$$\begin{array}{lllll} x_1=\frac{2^{m-l}-1}{2^l-1}, &x_2=\frac{2^{l-1}(2^{m-l}-1)}{2^l-1},
&x_3=\frac{q-q'+\delta r-\delta' r'}{4(q'+\delta'r'+1)},
&x_4=\frac{q-q'-\delta r+\delta' r'}{4(q'-\delta'r'+1)},  &x_5=\frac{2^{m-l}(2^{m-l}-1)}{2^l(2^l-1)}.\end{array}$$
\end{lem}
\demo
Since $K_1$ fixes $|\N_{M_0}(K_1):K_1|=|Q_0.K_0:Q_0|=2^l-1$ points in each suborbit whose point stabilizer is $M_0$-conjugate to $K_1$, we obtain $x_1=\frac{\Fix(K_1)-1}{|\N_{M_0}(K_1):K_1|}=\frac{2^{m-l}-1}{2^l-1}$.

Since $K_2$ fixes  $|\N_{M_0}(K_2):K_2|=2$ points in each suborbit whose point stabilizer is $M_0$-conjugate to $K_2$, we obtain $x_2=\frac{\Fix(K_2)-1}{2}=\frac{\frac{q-1}{q'-1}-1}{2}=\frac{q-q'}{2(q'-1)}=\frac{2^{l-1}(2^{m-l}-1)}{2^l-1}$.

Since $K_3$ fixes $|\N_{M_0}(K_3):K_3|=4$ points in each suborbit whose point stabilizer is $M_0$-conjugate to $K_3$, we obtain $x_3=\frac{|\Fix(K_3)-1|}{|\N_{M_0}(K_3):K_3|}=\frac{q-q'+\delta r-\delta' r'}{4(q'+\delta'r'+1)}$. Similarly, $x_4=\frac{q-q'-\delta r+\delta' r'}{4(q'-\delta'r'+1)}$.

Since $K_5=\Z(K_1)$ is the unique elementary abelian subgroup of $K_1$ of order $q'$, the group $K_5$ fixes $|\N_{M_0}(K_5):\N_{K_1}(K_5)|=|Q_0.K_0:Q_0|=2^l-1$ points in each suborbit whose point stabilizer is $M_0$-conjugate to $K_1$. Also, $K_5$
fixes $|\N_{M_0}(K_5):K_5|=|Q_0.K_0:Z(Q_0)|=2^l(2^l-1)$ points in each suborbit whose point stabilizer is $M_0$-conjugate to $K_5$. Therefore,  $x_5=\frac{\Fix(K_5)-1-x_1\times (2^l-1)}{2^l(2^l-1)}=\frac{2^{m-l}(2^{m-l}-1)}{2^l(2^l-1)}$.
\qed
\vskip 3mm
As a consequence, we have the following corollary.
\begin{cor}  \label{regular-Su}  Let    $\G_0$ be  the union of all regular suborbits of $T$ on $[T:M_0]$. Then
$  |\G_0|=|\O_0|-1-\sum_{i=1}^5\frac{|M_0|}{|K_i|}x_i.$
Moreover, as will be shown in the next subsection, $|\G_0|\gneqq \frac {|\O_0|}2$.
\end{cor}

\subsection{The proof of Theorem~\ref{main} for $\soc(G)=\Sz(q)$}
 Recall that $q=2^m$ and $q'=2^l$, where $m,l\geq3$ are odd and $e := m/l$ is a prime; $T=\Sz(q)$, $G=T\rtimes\lg f_1\rg,$ where $f_1\in\lg f\rg$ and $|f_1|=m_1$, for some $m_1\mid m$; $M_0=\Sz(q')$ and $M=M_0.\lg f_1\rg.$
In this subsection, we shall consider the action of $G$ on $\O=[G:M]$, while
  the action of $T$ on $[T:M_0]$ is identified with the action of $T$  on $[G:M]$, with a point stabilizer $M_0.$ Fix $\a=M\in\Omega$.
   Let $\G$ and $\G_0$ be the unions of all regular suborbits of $G$ and $T$,  respectively, relative to $\a$.  Then the size of $\G_0$  has been expressed  in Corollary~\ref{regular-Su}.

Each regular suborbit of $G$ is a disjoint union of $m_1$ regular suborbits of $T$. However, some regular suborbits of $T$ may lie in nonregular suborbits of $G$. To describe these, let $\Delta$ be a regular suborbit of $T$, choose $\b\in\Delta$, and let $M_{\Delta}$ be its setwise  stabilizer in $M$. Since $M_0$ acts regularly on $\Delta$, we have $M_{\Delta}=M_0M_{\beta},$ $M_0\cap M_{\beta}=1$. Thus, $M_{\beta}\cong M_{\Delta}/M_0$. If $M_{\beta}\neq 1$, it contains a subgroup $\lg vf'\rg$ of prime order $\ell$, where $v\in M_0$ and $f'\in\lg f_1\rg$ has order $\ell$. Consequently, $\Delta$ lies in a nonregular  suborbit of $G$ if and only if it is setwise stabilized by some subgroup $L\leq M$ of prime order with $L\nleq M_0$.

  We now introduce  two parameters.

\vskip 3mm
{\it Let $n'(\ell)$ be the number of regular $T$-suborbits  which are setwise stabilized by a subgroup $L\leq M$ of prime order $\ell$ with $L\nleq M_0$; and let $n'$ be the total number of distinct regular $T$-suborbits setwise stabilized by some subgroup $L\leq M$ of prime order with $L\nleq M_0$.}

\vskip 3mm
 Clearly, $|\G|=|\G_0|- n'|M_0|$. Therefore, Theorem~\ref{main} holds if $|\G|-\frac{|\O|}{2}\gneqq 0$, that is  \vspace{-5pt}$$|\G_0|-n'|M_0|-\frac{|\O|}{2}\gneqq 0.$$
In the next two lemmas, we give upper bounds for $n'(\ell)$ and $n'$, respectively.

\begin{lem}\label{e'=e} With the above notation, let  $\ell$ be a prime divisor of $m_
1$ and let  $L=\lg vf'\rg\cong\ZZ_\ell $, where $v\in M_0$ and $ f'\in\lg f\rg$. Suppose that $L$ fixes a point in a regular $T$-suborbit. Then the following bounds hold.
\begin{enumerate}
\item[\rm(1)] $n'(\ell)\leq \frac{2^{\frac{2m}\ell}(2^{\frac{2m}\ell}+1)(2^{\frac{m}\ell}-1)}{[2^{\frac{2m}{\ell e}}(2^{\frac{2m}{\ell e}}+1)(2^{\frac{m}{\ell e}}-1)]^2},$  if $\ell\neq e$;
\item[\rm(2)]  $n'(\ell)\leq \frac{65}{224}(\ell-1)^2q'^3$ if $\ell=e$.
\end{enumerate}
\end{lem}
\demo Let $e\mid m$,  $m_1\mid m$ and $\ell\mid m_1$, where $e$ and $\ell$ are primes. Let $L=\lg vf'\rg\cong \ZZ_\ell$, which fixes a point in a regular $T$-suborbit.  Recalling $M_0=\Sz(2^l)$ and $L\le M=M_0.\lg f_1\rg $ where $|f_1|=m_1$,   we consider the cases $\ell\ne e$ and $\ell=e$  separately.
\vskip 3mm

(1) Suppose that $\ell\neq e$. Then $\ell e\mid m$. We first compute $\Fix(L)$. In the notation of Proposition~\ref{fieldconjugate},  we have  $\Inndiag (M_0)=\Inn(M_0)=M_0$. By this proposition,    $vf'$ is conjugate to $f'$ by an element of $\Inndiag (M_0)=M_0$. Therefore,  $L=\lg vf'\rg $ is $M$-conjugate to $\lg f'\rg $.
  Now, by Proposition \ref{man}, we have
  \vspace{-5pt}$$\begin{array}{lll}\Fix(L)&=&\Fix(\lg f'\rg )=|\N_G(\lg f' \rg)|/|\N_M(\lg f'\rg)|=|\Sz(2^{\frac{m}{\ell}})\lg f_1\rg|/|\Sz(2^{\frac{m}{\ell e}})\lg f_1\rg|\\
  &=&\frac{2^{\frac{2m}\ell}(2^{\frac{2m}\ell}+1)(2^{\frac{m}\ell}-1)}{2^{\frac{2m}{\ell e}}(2^{\frac{2m}{\ell e}}+1)(2^{\frac{m}{\ell e}}-1)}.\end{array}$$
\f Let $\a^{sM_0}$ be a regular $T$-suborbit containing a point fixed by $\lg vf'\rg $. Consider the action of $M_0\lg f'\rg $ on $\a^{sM_0}$.  Since $L=\lg vf'\rg $ is conjugate to $\lg f'\rg $ in $M_0\lg f'\rg$, we know  that  $\lg f'\rg $ fixes some points on $\a^{sM_0}$.  In fact,  $\lg f'\rg $ fixes exactly $|\N_{M_0\lg f'\rg}(\lg f'\rg):\lg f'\rg|={2^{\frac{2m}{\ell e}}(2^{\frac{2m}{\ell e}}+1)(2^{\frac{m}{\ell e}}-1)}$ points. Therefore,
\vspace{-5pt}$$n'(\ell)\leq \frac{\Fix(\lg f'\rg )}{2^{\frac{2m}{\ell e}}(2^{\frac{2m}{\ell e}}+1)(2^{\frac{m}{\ell e}}-1)}=   \frac{2^{\frac{2m}\ell}(2^{\frac{2m}\ell}+1)(2^{\frac{m}\ell}-1)}{[2^{\frac{2m}{\ell e}}(2^{\frac{2m}{\ell e}}+1)(2^{\frac{m}{\ell e}}-1)]^2}.$$

\vskip 3mm (2) Suppose that $\ell=e$.  Then $e\di m_1$  and   $[M_0, f']=1$ and we need to deal with three cases, separately: $\ell\nmid|M_0|$, $\ell\mid q'\pm\delta'r'+1$ and $\ell\mid q'-1$, where $q'=2^{\frac m\ell}$ and $\d'=\pm 1$ such that $5\di (q'+\d' r'+1)$.

\vskip 3mm

(2.1)  $\ell\nmid |M_0|$

In this case, $L=\lg f'\rg$. Since $\N_G(\lg f'\rg)=\N_M(\lg f'\rg)=M_0\lg f_1 \rg$, we get $\Fix(L)=1$ and so
 $n'(\ell)=0$.

\vskip 3mm

(2.2)  $\ell\mid q'\pm\delta'r'+1$

 Let $\lg u_0\rg $ be a subgroup of order $\ell$ in $M_0$. Up to $M_0$-conjugacy, we may assume $L\le \lg u_0\rg \times \lg f'\rg \cong \ZZ_\ell^2$. Then
 $L$ is $M_0$-conjugate to $L_i:=\lg u_0^if'\rg$, where $0\leq i \leq \ell-1$. By
Proposition~\ref{fieldconjugate} again, $\lg u_0^if'\rg$  is  conjugate to $\lg f'\rg$ in $G$, and so $|\N_G(L_0)|=|\N_G(L_i)|=|M|$ $(0\leq i\leq \ell-1)$. Since $\N_{M_0}(\lg u_0\rg)\cong  \ZZ_{q'\pm\delta'r'+1}\rtimes\ZZ_4$, and $\C_{M_0}( u_0)\cong \ZZ_{q'\pm\delta'r'+1}$, there are $\frac{\ell-1}{4}$ $M_0$-conjugacy classes of elements of order $\ell$.
Some of these classes may merge under conjugation by $M$, as $M_0\unlhd M.$ Let $J_1,\cdots,J_t\in M_0$ be representatives for the $M$-conjugacy classes of elements of order $\ell$, and suppose that the $M$-conjugacy class of $J_j$ contains $a_j$ of the above $M_0$-classes. Therefore, $\sum\limits_{j=1}^ta_j=\frac{\ell-1}{4}$.
Since $\lg u_0^if'\rg^g=\lg (u_0^i)^gf'\rg$ for any $g\in M,$ we have $\N_M(\lg u_0^if'\rg)=\C_M(u_0^i)$
and the $M$-conjugacy of these subgroups corresponds precisely to $M$-conjugacy of the elements $u_0^i$.
Note that $\frac{|M|}{|\C_M(J_j)|}=a_j\frac{|M_0|}{|\C_{M_0}(J_j)|}$.
By Proposition~\ref{man}, we have
\vspace{-5pt}$$\begin{array}{lll}&&\Fix(\lg u_0^if'\rg )=\frac{|\N_G(\lg f'\rg)|}{|\N_M(\lg f'\rg)|}+\sum_{j=1}^{t}\frac{|M|}{|\C_M(J_j )|}
=1+\sum_{j=1}^{t}a_j\frac{|M_0|}{|\C_{M_0}(J_j )|}\\
 &&= \frac{|M|}{|M|}+\frac{1}{4}(\ell-1)\frac{|M_0|}{|\C_{M_0}(\lg u_0\rg)|}=1+\frac{1}{4}(\ell-1)\frac{q'^2(q'^2+1)(q'-1)}{q'\pm\delta'r'+1} .\end{array}$$
Since $|\N_{M_0\lg f'\rg}(\lg f'\rg):\lg f' \rg|=|M_0|>\Fix(L_0)$, we have that $\lg f' \rg$ can not fix a point in any regular $T$-suborbit.
Now we will consider the case $\lg u_0^if' \rg$ $(i\neq 0)$. Except for $\a$,
$\lg u_0^i f' \rg$ fixes at most $\frac{1}{4}(\ell-1)\frac{q'^2(q'^2+1)(q'-1)}{q'\pm\delta'r'+1} $ points. In each regular $T$-suborbit containing a point fixed by $L_i$, $\lg u_0^i f' \rg$ fixes $|\N_{M_0\lg f'\rg}(\lg u_0^if' \rg):\lg u_0^if' \rg|\geq (q'\pm\delta'r'+1)$ points and so
$n'(\ell)\leq\frac{1}{16}(\ell-1)^2\frac{q'^2(q'^2+1)(q'-1)}{(q'\pm\delta'r'+1)^2}\leq \frac{65}{256}(\ell-1)^2q'^2(q'-1),$ as $l\geq3$ and $r'=2^{(l+1)/2}.$

\vskip 3mm
(2.3) $\ell\mid q'-1$

 Let $\lg u_0\rg $ be a subgroup of order $\ell$ in $M_0$. Up to $M_0$-conjugacy, we may assume that $L\le \lg u_0\rg \times \lg f'\rg \cong \ZZ_\ell^2$. Then
 $L$ is conjugate to $L_i:=\lg u_0^if'\rg$, where $0\leq i \leq \ell-1$. By
Proposition~\ref{fieldconjugate} again,    $\lg u_0^if'\rg$  is  conjugate to $\lg f'\rg$ in $G$, and so $|\N_G(L_0)|=|\N_G(L_i)|=|M|$ $(0\leq i\leq \ell-1)$.  Since $\N_{M_0}(\lg u_0\rg)\cong  \ZZ_{q'-1}\rtimes\ZZ_2$, and $\C_{M_0}( u_0)\cong \ZZ_{q'-1}$, there are $\frac{\ell-1}{2}$ $M_0$-conjugacy classes of elements of order $\ell$.
Apply the $M$-conjugacy-class calculation in (2.2). If the $M$-class represented by $J_j$ contains $a_j$ of these $M_0$-classes, then $\sum\limits_{j}a_j=\frac{\ell-1}{2}$
and
$\frac{|M|}{|\C_M(J_j)|}=a_j\frac{|M_0|}{|\C_{M_0}(J_j)|}$.
 Now, by Proposition~\ref{man}, we get
\vspace{-5pt}$$\begin{array}{lll}&&\Fix(\lg u_0^if'\rg )=\frac{|\N_G(L_0)|}{|\N_M(L_0)|}+\sum_{j}\frac{|M|}{|\C_M(J_j)|}=1+\sum_{j}a_j\frac{|M_0|}{|\C_{M_0}(J_j)|}\\
  &&=\frac{|M|}{|M|}+\frac{1}{2}(\ell-1)\frac{|M_0|}{|\C_{M_0}(\lg u_0\rg)|}=1+\frac{1}{2}(\ell-1)\frac{q'^2(q'^2+1)(q'-1)}{q'-1} .\end{array}$$
Since $|\N_{M_0\lg f'\rg}(\lg f'\rg):\lg f' \rg|=|M_0|>\Fix(L_0)$, we have that $\lg f' \rg$ can not fix a point in any regular $T$-suborbit.

Now we will consider the case $\lg u_0^if' \rg$ $(i\neq 0)$. Except for $\a$,
$\lg u_0^i f' \rg$ fixes at most $\frac{1}{2}(\ell-1)\frac{q'^2(q'^2+1)(q'-1)}{q'-1} $ points. In each regular $T$-suborbit containing a point fixed by $L_i$,  $\lg u_0^i f' \rg$ fixes $|\N_{M_0\lg f'\rg}(\lg u_0^i f' \rg):\lg u_0^i f' \rg|\geq (q'-1)$ points and so
$n'(\ell)\leq\frac{1}{4}(\ell-1)^2\frac{q'^2(q'^2+1)(q'-1)}{(q'-1)^2}=\frac{1}{4}(\ell-1)^2q'^3\times\frac{(1+q'^{-2})}{1-q'^{-1}}\leq\frac{65}{224}(\ell-1)^2q'^3 $, as $q'\geq8$.

Taking the maximal $n'(\ell)$ in the above three cases, we get the result as desired.
\qed
\vskip 3mm
To facilitate the following proof, we write $l=e_1^{o_1}e_2^{o_2}\cdots e_n^{o_n}\ge 3$, where $3\leq e_1<e_2<\cdots<e_n$ are primes and $o_i$ are positive integers.
\begin{lem}\label{n'} With the above notation,  we have $$n'\leq \frac{65}{56}q'^e+\frac{9}{8}q'^{\frac{5}{3}e-\frac{7}{3}}.$$
\end{lem}
\demo Since $l\geq 3^{n-1}e_1$, we have $n\leq 1+\log_3\frac{l}{e_1}\leq 1+\log_3l$. Hence,
\vspace{-5pt}$$n'\leq n'(e)+\sum_{i=1}^nn'(e_i)~~for~~e\nmid l, \quad n'\leq n'(e)+\sum_{e_i\neq e}n'(e_i)~~for ~~e\mid l.$$
Moreover, by Lemma \ref{e'=e} we have
\vspace{-5pt}$$\begin{array}{ll}
n'(e)\leq  \frac{65}{224}(e-1)^2q'^3,&n'(e_i)\leq\frac{2^{2l_ie}(2^{2l_ie}+1)(2^{l_ie}-1)}{[2^{2l_i}(2^{2l_i}+1)(2^{l_i}-1)]^2}\leq\frac{9}{2}\times2^{5l_1(e-2)},
\end{array}$$
 where $l=l_ie_i$ with $1\leq i\leq n$ and $e_i\neq e$. Since $m=ee_1^{o_1}e_2^{o_2}\cdots e_n^{o_n}$, we get that
  \vspace{-5pt}$$\begin{array}{lll}
  n'&\leq& \frac{65}{224}(e-1)^2q'^3+\frac{9n}{2}\times2^{5l_1(e-2)}\leq \frac{65}{224}(e-1)^2q'^3+\frac{9(1+\log_3l)}{2}\times2^{5l_1(e-2)}\\
  &\leq&\frac{65}{56}q'^e+\frac{9}{8}q'^{\frac{5}{3}e-\frac{7}{3}}\quad(\text{as}\,\,\frac{(e-1)^2}{q'^e}\leq\frac{4}{q'^3},\,\,\frac{1+\log_3l}{2^l}\leq \frac{1}{4}\,\,\text{and}\,\,l_1\leq\frac{l}{3}).
  \end{array}$$
 \qed

Now we are ready to prove  Theorem~\ref{main} for  $\soc(G)=\Sz(q)$.
\begin{theorem}\label{main-Sz}
Let $G=\Sz(q):\lg f_1\rg $ and $M=\Sz(q'):\lg f_1\rg $,  where $q=2^m$,  $q'=2^l$, $l\di m$,  $|f_1|=m_1\di m$ with $m,l\geq 3$ odd integers.    Consider the primitive action of $G$ on $\O:=[G:M]$ by right multiplication. Then $b(G)=2$ and the BG-Conjecture holds.
\end{theorem}
\demo By the preceding arguments, to prove the theorem, it suffices to show  $A:=|\G_0|-n'|M_0|-\frac{|\O|}{2}\gneqq 0.$
From $|\Omega_0| = |\Omega|$, Corollary~\ref{regular-Su}, and Lemma~\ref{n'}, it follows that
\vspace{-5pt}$$A=
 \frac{|\O|}{2}-1-\sum_{i=1}^5\frac{|M_0|}{|K_i|}x_i-n'|M_0|.$$
Put $D:=|\G_0|-\frac{|\Omega|}{2}=\frac{|\O|}{2}-1-\sum_{i=1}^5\frac{|M_0|}{|K_i|}x_i$.
Then $A=D-n'|M_0|$.
Substituting the explicit formulas for $|\Omega|$, $|M_0|$, $|K_i|$, and $x_i$ (from the preceding lemmas), we obtain the following lower bound for $D$:
\vspace{-5pt}$$\begin{array}{lll}
 D&=& \frac{1}{2}\frac{2^{2m}(2^{2m}+1)(2^m-1)}{2^{2l}(2^{2l}+1)(2^l-1)}-1-\frac{|M_0|}{|K_1|}x_1-\frac{|M_0|}{|K_2|}x_2-\frac{|M_0|}{|K_3|}x_3-\frac{|M_0|}{|K_4|}x_4-\frac{|M_0|}{|K_5|}x_5\\
&=&\frac{1}{2}\frac{2^{2m}(2^{2m}+1)(2^m-1)}{2^{2l}(2^{2l}+1)(2^l-1)}-1-\frac{2^{2l}(2^{2l}+1)(2^l-1)}{2^{2l}}\frac{2^{m-l}-1}{2^l-1}-\frac{2^{2l}(2^{2l}+1)(2^l-1)}{2^l-1}\frac{2^{l-1}(2^{m-l}-1)}{2^l-1}\\
&&-\frac{2^{2l}(2^{2l}+1)(2^l-1)}{q'+\delta'r'+1}\frac{q-q'+\delta r-\delta' r'}{4(q'+\delta'r'+1)}
-\frac{2^{2l}(2^{2l}+1)(2^l-1)}{q'-\delta'r'+1}\frac{q-q'-\delta r+\delta' r'}{4(q'-\delta'r'+1)}
-\frac{2^{2l}(2^{2l}+1)(2^l-1)}{q'}\frac{2^{m-l}(2^{m-l}-1)}{2^l(2^l-1)}\\
&\ge&\frac{511}{1040}q'^{5e-5}-1-\frac{65}{64}q'^{e+2}-\frac{65}{112}q'^{e+3}-2.1q'^{e+3}-\frac{65}{64}q'^{2e}\quad(\text{as}\,\, e,l\geq3)\\
&\geq&\frac{511}{1040}q'^{5e-5}-2.681q'^{e+3}-\frac{65}{64}q'^{e+2}-\frac{65}{64}q'^{2e}-1>0.
\end{array}$$
Using the upper bound $|M_0| \le \frac{65}{64} q'^5$ together with the lower bound for $D (>0)$ obtained above and the upper bound for $n'$ from Lemma~\ref{n'}, we derive a lower bound for $A/|M_0|$:
\vspace{-5pt}$$\begin{array}{lll}
\frac{A}{|M_0|}&=&\frac{D}{|M_0|}-n'\geq B,
\end{array}$$
where $B:=0.48q'^{5e-10}-2.64q'^{e-2}-q'^{e-3}-q'^{2e-5}-\frac{64}{65}q'^{-5}-(\frac{65}{56}q'^e+\frac{9}{8}q'^{\frac{5}{3}e-\frac{7}{3}})$.
Since $|M_0| > 0$, proving $A > 0$ is equivalent to proving $A/|M_0| > 0$. Consequently, it suffices to show that the lower bound $B$ is positive.
Since $e\geq3$ and $q'\geq8$, we have
\vspace{-5pt}$$\begin{array}{lll}
\frac{B}{q'^{2e}}&=&0.48q'^{3e-10}-2.64q'^{-e-2}-q'^{-e-3}-q'^{-5}-\frac{64}{65}q'^{-2e-5}-\frac{65}{56}q'^{-e}-\frac{9}{8}q'^{-\frac{1}{3}e-\frac{7}{3}}\gneqq0,
\end{array}$$
as desired.
\qed

\vskip 3mm
\begin{center}{\large\bf Acknowledgements}\end{center}
\vskip 2mm
The authors acknowledge support from the National Natural Science Foundation of China: H. Chen (12301446, 12571362),
S. Du (12471332).

\end{document}